\documentclass[11pt]{article}

\usepackage{amssymb}
\usepackage{graphicx}
\usepackage{color}
\usepackage{amsmath,amsthm,amsfonts,epsfig,setspace}

\newtheorem{thm}{Theorem}[section]
\newtheorem{rmk}{Remark}[section]

\newtheorem{definition}{Definition}[section]
\newtheorem{lem}{Lemma}[section]

\newtheorem{prop}{Propsition}

 \def\p{\partial}
\def \Vh0{\stackrel{\circ}{V}_h} \def\to{\rightarrow}
   
\def\Om{\Omega}  
\newcommand{\q}{\quad}   
\def\l{\label}  \def\f{\frac}  
\def\D{\end{document}}

\def\l|{\left|}
\def\r|{\right|}

\newcommand{\R}{\mathbb{R}}

\newcommand{\lc}
{\mathrel{\raise2pt\hbox{${\mathop<\limits_{\raise1pt\hbox
{\mbox{$\sim$}}}}$}}}

\newcommand{\gc}
{\mathrel{\raise2pt\hbox{${\mathop>\limits_{\raise1pt\hbox{\mbox{$\sim$}}}}$}}}

\newcommand{\ec}
{\mathrel{\raise2pt\hbox{${\mathop=\limits_{\raise1pt\hbox{\mbox{$\sim$}}}}$}}}

\def\be{\begin{equation}} \def\ee{\end{equation}}

\def\bea{\begin{eqnarray}}  \def\eea{\end{eqnarray}}

\def\beas{\begin{eqnarray*}} \def\eeas{\end{eqnarray*}}

\def\bn{\begin{enumerate}} \def\en{\end{enumerate}}

\def\bd{\begin{description}} \def\ed{\end{description}}

\title{Sensitivity Analysis of an Inverse Problem for the Wave Equation with Caustics}

\author{Gang Bao\thanks{Department of Mathematics, Zhejiang University, Hangzhou 310027,
China; Department of Mathematics, Michigan
State University, East Lansing, MI 48824. The research was supported
in part by the NSF grants DMS-0908325, DMS-0968360,
DMS-1211292,  the ONR grant N00014-12-1-0319, a Key
Project of the Major Research Plan of NSFC (No. 91130004), and a
special research grant from Zhejiang University. Email: {\bf bao@math.msu.edu}}
\and Hai Zhang\thanks{Department of Mathematics, Michigan State University, East Lansing, MI 48824. Email: {\bf zhangh20@msu.edu}}
}

\begin{document}
\maketitle

\begin{abstract}
The paper investigates the sensitivity of the inverse problem of recovering the velocity field in a bounded domain from the boundary dynamic Dirichlet-to-Neumann map (DDtN) for the wave equation. Three main results are obtained: (1) assuming that
two velocity fields are non-trapping and are equal to a constant near the boundary, it is shown that
the two induced scattering relations must be identical if their corresponding DDtN maps are sufficiently close;
(2) a geodesic X-ray transform operator with matrix-valued weight is introduced
by linearizing the operator which associates each velocity field with its induced Hamiltonian flow.
A selected set of geodesics whose conormal bundle can cover the cotangent space at an interior point is
used to recover the singularity of the X-ray transformed function at the point;  a local stability estimate is established for this case.
Although fold caustics are allowed along these geodesics,
it is required that these caustics contribute to a smoother term in the transform than the point itself.
The existence of such a set of geodesics is guaranteed under some natural assumptions in dimension greater than or equal to three
by the classification result on caustics and regularity theory of Fourier Integral Operators.
The interior point with the above required set of geodesics is called ``fold-regular'';
(3) assuming that a background velocity field with every interior point fold-regular is fixed and
another velocity field is sufficiently close to it and satisfies a certain orthogonality condition,
it is shown that if the two corresponding
DDtN maps are sufficiently close then they
must be equal.
\end{abstract}

\medskip


\section{Introduction}
This paper is concerned with the sensitivity (or stability) of the inverse problem of recovering the velocity field in a domain from the boundary dynamic Dirichlet-to-Neumann map (DDtN) in the wave equation. Let $\Omega$ be a bounded strictly convex smooth domain in $\R^d$, $d\geq 2$, with boundary $\Gamma$. Let $c(x)$ be a velocity field in $\Omega$ which characterizes
the wave speed in the medium and let $T$ be a sufficiently large positive number.
We consider the following wave equation system:
\begin{eqnarray}
\frac{1}{c^2(x)}u_{tt}-\Delta u&=&0,   \quad (x,t)\in \R^d \times(0, T)  \label{wave1}\\
u(0,x)=u_t(0,x)&=&0,                   \quad  x\in \Omega, \label{wave2}\\
u(x,t)&=&f(x,t),                   \quad   (x,t)\in \Gamma \times (0,T). \label{wave3}
\end{eqnarray}
For each $f\in H_0^1([0, T]\times \Gamma)$, it is known that (see for instance \cite{LLT86}) there exists a unique solution $u\in C^1(0, T; L^2(\Om))\bigcap C(0, T; H^1(\Om))$, and furthermore
$\frac{\partial u}{\partial \nu}\in L^2([0, T]\times \Gamma)$, where $\nu$ is the unit outward normal to the boundary.
The DDtN map $\Lambda_c$ is defined by
$$
 \Lambda_c(f):= \frac{\partial u}{\partial \nu}|_{[0, T]\times \Gamma}.
$$
The inverse problem is to recover the velocity function $c$ from the DDtN map $\Lambda_c$. The uniqueness of the inverse problem is solved by the boundary control method first introduced by Belishev in \cite{B87}. The method can also be used to solve the uniqueness for more
general problems, for instance, the anisotropic medium case. We refer to \cite{B97}, \cite{BK92}, \cite{B07}, \cite{KKL01}, \cite{R87} and the references therein for more discussions.

We are interested in the sensitivity question for the above inverse problem. Namely, we want to investigate how sensitive or stable is it to recover the velocity field from the DDtN map and characterize how a small change in the DDtN map affects the recovered velocity field.

The above inverse problem can be viewed as a special case of the problem of recovering a Riemannian metric on a Riemannian manifold. Indeed, it corresponds to the case when the metrics are restricted to the class of those which are conformal to the Euclidean one.
The inverse problem of recovering a Riemannian metric
has been extensively studied in the literature. The uniqueness is proved by Belishev and Kurylev in \cite{BK92} by using the boundary control method. However, as pointed out in \cite{SU05IMRN}, their approach is unlikely to give a stability estimate since it uses in an essential way a unique continuation property of the wave equation.

The first stability result on the determination of the metric from the DDtN map was given by Stefanov and Uhlmann in \cite{SU98}, where they proved conditional stability of H$\ddot{o}$lder type for metrics close enough to the Euclidean one in $C^k$ for $k\gg 1$ in three dimensions. Later, they extended the stability result to generic simple metrics, \cite{SU05IMRN}. An important feature of their approach is to first derive a stability estimate of recovering the boundary distance function from the DDtN map and then apply existing results from the boundary rigidity problem in geometry. Their approach was extended by Montalto in \cite{M-arxiv} to study the more general problem of determine a metric, a co-vector and a potential simultaneously from the DDtN map, and a similar  H$\ddot{o}$lder type conditional stability result was obtained. The stability of the inverse problem of determining the conformal factor to a fixed simple metric was studied by Bellassoued and Ferreira in \cite{BF11}. They proved the H$\ddot{o}$lder type conditional stability result for the case when the conformal factors are close to one. We comment that the result in \cite{BF11} holds for all simple metrics. For other stability results on the related problems, we refer to the references in \cite{M-arxiv}.

We emphasize that all of the above stability results deal with the case when the metrics are simple. To our best knowledge, no stability result is available in the general case when the metrics are not simple. This paper is devoted to the study of the general case when the metric induced by the velocity field is not simple. To avoid technical complications due to the boundary, we restrict our study to situation when the velocity fields are equal to one near the boundary. From this point of view, our results can be regarded as interior estimates. We refer to \cite{SU05IMRN}, \cite{SU08AJM} and the references therein for useful boundary estimates.

We now give a brief account of the approach and results in the paper.
We first derive a sensitivity result of recovering the scattering relation from the DDtN map.
Our result shows that two scattering relations must be identical if the two corresponding DDtN maps are sufficiently close in some suitable norm. Equivalently, any arbitrarily small change in the scattering relation can imply a certain change in the DDtN map.  To our best knowledge, this seems to be the first sensitivity result for the inverse problem in the non-simple metric case. Moreover,
our result is fundamentally different from those in the literature where Lipschitz, H$\ddot{o}$lder or logarithmic estimates are derived, see for examples \cite{BY09}, \cite{AS90}, \cite{S90}, \cite{SU98}, \cite{SU04Duke}, \cite{SU05JAMS}, \cite{SU05IMRN}, \cite{M-arxiv} and \cite{Isakov06}. This is the reason the term ``sensitivity analysis'' is used instead of ``stability estimates'' throughout the paper.
We remark that when the geometry induced by the velocity field is simple, the scattering relation is equivalent to the boundary distance function.
In that case, a H$\ddot{o}$lder type interior stability estimate for recovering the boundary distance function from the DDtN map has been established in
\cite{SU05IMRN}. Compared to the H$\ddot{o}$lder type result, our result is much stronger. As a consequence, a stronger result than H$\ddot{o}$lder type stability may be proved for the results in \cite{SU05IMRN}, \cite{BF11} and \cite{M-arxiv}.
Our approach is based on Gaussian beam solutions to the wave equation, which are capable of dealing with caustics, major obstacles
to the construction of classic geometric-optics solutions.
We refer to \cite{R83} and \cite{KKL01} for more discussions on Gaussian beams and its applications.

Once the above sensitivity result of recovering the scattering relation from its associated DDtN map is established, the sensitivity of the inverse problem of recovering the velocity field from its induced DDtN map is reduced to the uniqueness issue of recovering the velocity field from its induced scattering relation. This is a special case of the lens rigidity problem in geometry when the metrics are restricted to a conformal class (see \cite{Mi81} and \cite{SU09JDG}).
With regarding to the lens rigidity problem, Stefanov and Uhlmann (\cite{SU09JDG}) proved uniqueness up to diffeomorphisms fixing the boundary for metrics a priori close to a generic ``regular'' one. There ``regular''  means that there is a set of geodesics without conjugate points whose conormal bundles can cover the contangent bundle of the underlying manifold. Their approach is based on investigating the X-ray transform which is obtained from linearizing the boundary distance function with respect to the metric. The lens rigidity problem remains open for more general non-simple metrics.

In this paper, motivated by \cite{SUAPDE} which studies the geodesic X-ray transform with fold caustics, we introduce the ``fold-regular'' class (see Subsection \ref{subsec-thm3} for definition) of velocity fields (or equivalently conformal metrics) which generalizes the above ``regular'' class,  and study an linearized problem of a variant of the lens rigidity problem when the metrics are restricted to those conformal to the Euclidean one. Our approach is based on linearizing the Hamiltonian flow with respect to the velocity field. This gives us the advantage of handling geodesics with caustics in comparison to the approach by using the
boundary distance function. More specifically,
we linearize the operator which maps $c$ to $\mathcal{H}^T_{c}|_{S^*\R^d}$ and obtain a geodesic X-ray transform operator $\mathfrak{I}_c$ with matrix-valued weight.
We study the inverse problem of recovering a vector-valued function $f$ from its X-ray transform $\mathfrak{I}_cf$. For a fixed interior point $x$,
we use a carefully selected set of geodesics whose conormal bundle can cover the cotangent space $T_x^*\R^d$ to recover the singularity of $f$
at $x$.
We allow fold caustics along these geodesics, but require that these caustics contribute to a smoother term in the transform than $x$ itself. It is still an open problem to show that such a set of geodesics exists generically for a general velocity field with caustics. But we draw evidence from the classification result on caustics and regularity theory of Fourier Integral Operators (FIOs) to show that it is the case under some natural assumptions in the dimensions equal or greater than three.
We call the interior point with the above set of geodesics ``fold-regular''. See Section \ref{subsec-discussion} for discussions on the new concept ``fold-regular''.  A local stability estimate is derived near a fold-regular point.
We think that this local stability estimate and the approach presented here may hold the key for analyzing the stability of the lens rigidity problem with general non-simple metrics, which is completely open at present. Initial progress has been made along this direction in \cite{thesis-zhang}.


Finally, we combine the stability result on the X-ray transform and the sensitivity result on recovering
the scattering relation from the DDtN map to obtain a sensitivity result for the inverse problem.

The paper is organized as follows. In Section 2, we present some preliminaries. In Section 3, we introduce the main results in the paper.
Section 4 is devoted to the construction of Gaussian beam solutions to the wave equation.
The Gaussian beam solutions are used in Section 5 to prove the sensitivity result of determining scattering relation from DDtN map.
In Section 6, we discuss the concept ``fold-regular'' and prove a local stability estimate for the geodesic X-ray transform $\mathfrak{I}_c$.
In Section 7, we prove the sensitivity result of recovering the velocity field from the DDtN map.

\medskip

Throughout, we use the following conventions:
\bn
\item
Let $f$ and $g$ be two elements in a Hilbert space, then $\langle f, g \rangle$ stands for their inner product;

\item
Let $M_1$ and $M_2$ be two matrices (including vectors which can be regarded as single column or single row matrices), then the product of $M_1$ and $M_2$ is denoted by $M_1\cdot M_2$. Sometimes, the dot is omitted for simplicity;

\item
Let $M$ be a matrix, then $M^{\dag}$ stands for its transpose. The same applies when $M$ is a linear operator. If $M$ is real and symmetric and $C$ a real number, then $M\geq C$ means that the
matrix $M-C\cdot Id$ is symmetric and positive definite. If $M$ is a complex matrix, then we use $\Re{M}$ for its real part and $\Im{M}$ for its imaginary part;

\item
Let $U$ and $V$ be two open set in a metric space, then $U\Subset V$ means that the closure of $U$, denoted by $\bar{U}$ is compact and is a subset of $V$;

\item
Let $C_1$ and $C_2$ be two positive numbers, then $C_1 \lesssim C_2$ means that $C_1\leq C\cdot C_2$ for some constant $C>0$ independent of $C_1$ and $C_2$.
\en

\section{Preliminaries}
In this section, we introduce some notations and definitions.
Let $\Omega$ be a strictly convex smooth domain in $\R^d$ with boundary $\Gamma$. Let $c$ be a smooth velocity field defined in $\Omega$ which is equal to one near the boundary. Then $c$ has natural extension to $\R^d$.
Throughout the paper, we always use the natural coordinate system of the cotangent bundle $T^*\R^d$ in which we write $(x, \xi)$ for the co-vector $\xi_j dx^j$ in $T^*_x\R^d$. For ease of notation, we also use $\xi$ for the co-vector $\xi_j dx^j$. The meaning of $\xi$ should be clear from the context. The velocity field $c$ introduces a Hamiltonian function $H_c(x,\xi)=\f{1}{2}c^2(x)(\xi_1^2+...\xi_d^2)$ to $T^*\R^d$.
It also defines a norm to each cotangent space $T^*_x\R^d$ by
$$
|\xi|=c(x)\sqrt{\xi_1^2+...\xi_d^2}, \q \mbox{for } ~\xi \in T^*_x\R^d.
$$


Denote the corresponding Hamiltonian flow by $\mathcal{H}_c^t$, i.e. for each $(x_0, \xi_0)\in T^*\R^d$,
$\mathcal{H}_c^t(x_0, \xi_0)= (x(t, x_0, \xi_0), \xi(t, x_0, \xi_0))$ solves the following equations:
\begin{eqnarray}
  \dot{x} &=& \f{\p H_c}{\p \xi}=c^2\cdot \xi, \quad x(0)=x_0 , \label{hamiltonian1}\\
  \dot{\xi} &=& -\f{\p H_c}{\p x}=-\f{1}{2}\nabla c^2 \cdot (\xi_1^2+...\xi_d^2), \quad \xi(0) = \xi_0. \label{hamiltonian2}
\end{eqnarray}
We call $(x(\cdot, x_0, \xi_0), \xi(\cdot, x_0, \xi_0))$ the bicharateristic curve emanating from $(x_0, \xi_0)$ and
$x(\cdot, x_0, \xi_0)$ the geodesic. By the assumptions on $c$, the flow $\mathcal{H}_c^t$ is defined for all $t\in \R$. Note that
the flow $\mathcal{H}_c^t$ is also well-defined on the cosphere bundle $S^*\R^d=\{(x, \xi): x\in \R^d, |\xi| =1\}$.

We say that a velocity field $c$ is non-trapping in $\Omega$ for time $T>0$ if the following condition is satisfied:
\be \label{nontrapping}
\mathcal{H}_c^T (S^*\Om) \bigcap S^*\Om = \emptyset.
\ee

Denote
\beas
S_{+}^* \Gamma&=& \{(x, \xi): x\in \Gamma, |\xi|= 1, \langle \xi, \nu(x) \rangle > 0 \}; \\
S_{-}^* \Gamma&=& \{(x, \xi): x\in \Gamma, |\xi|= 1, \langle\xi, \nu(x)\rangle < 0 \}.
\eeas

Assume that the velocity field $c$ is non-trapping in $\Omega$ for time $T$;
we now define the scattering relation  $\mathfrak{S}_c: S_{-}^*\Gamma \to S_{+}^*\Gamma$.
For each $(x_0, \xi_0)\in S_{-}^*\Gamma$, let $l(x_0, \xi_0)$ be the first moment that the geodesic
$x(\cdot, x_0, \xi_0)$ hits the boundary $\Gamma$. Define
$$
 \mathfrak{S}_c(x_0, \xi_0)= \mathcal{H}_c^{l(x_0, \xi_0)}(x_0, \xi_0).
$$

For future reference, we define $l_{-}:S^*\Omega \to (-\infty, 0]$ by letting $l_{-}(x, \xi)$ be the first negative moment
that the bicharacteristic curve $\mathcal{H}^t(x, \xi)$ hits the boundary
$S^*_{-}\Gamma$ and $\tau: S^*\Omega \to S^*_{-}\Gamma$ by
$$
 \tau(x, \xi)= \mathcal{H}^{l_{-}(x, \xi)}(x, \xi).
$$
We remark that $l_{-}(\cdot)$ and $\tau(\cdot)$ are well-defined by the assumption (\ref{nontrapping}).

We now introduce the class of admissible velocity fields that are considered in the paper.
\begin{definition}
Let $M_0$, $\epsilon_0$ and $T$ be positive numbers. A velocity field $c$ is said to belong to the
admissible class $\mathfrak{A}(M_0, \epsilon_0, \Om, T)$ if and only if the following three conditions are satisfied:
\begin{enumerate}
\item
$c\in C^3(\R^d)$, $0< \f{1}{M_0} \leq c \leq M_0$, and $\|c\|_{C^3(\R^d)} \leq M_0$;
\item
the support of $c-1$ is contained in the set $\Omega_{\epsilon_0}=: \{x\in \Omega: \mbox{dist}(x, \Gamma)> \epsilon_0\}$;
\item
the Hamiltonian $H_c$ is non-trapping in $\Omega$ for time $T$.
\end{enumerate}
\end{definition}

By Condition 2 above and the assumption that $\Omega$ is bounded and convex, it is easy to verify that any ray starting in the domain $\Omega_{\epsilon_0}$ intersect $\Gamma$ transversely. Moreover, there exist two small positive constants $\epsilon^*$ and $\epsilon_1$, both depending on $\epsilon_0$,
such that for any $(x_0, \xi_0) \in S^*_{-}\Gamma$, if
$\{\mathcal{H}^t(x_0, \xi_0): t\in (0, l(x_0, \xi_0))\} \bigcap S^*\Omega_{\epsilon_0} \neq \emptyset$, then
\bea
\langle \xi_0, \nu(x_0) \rangle &\leq & -\epsilon^*, \label{condition-ad1}\\
\langle \xi_1, \nu(x_1) \rangle &\geq & \epsilon^*, \label{condition-ad2}\\
 l(x_0, \xi_0) &\geq & \epsilon_1, \label{condition-ad3}
\eea
where $(x_1, \xi_1)= \mathfrak{S}_c(x_0, \xi_0)$.

\medskip

Finally, we remark that we set up the discussion in the paper in the cotangent space $T^*\R^d$.
But one can also set up the discussion in the tangent space $T\R^d$, see for instance \cite{S99}, \cite{SU09JDG}.
The equivalence of the two setups can be seen from the procedure of ``raising and lowing indices'' in Riemannian geometry.
We choose the cotangent setup mainly because the following three reasons. First, it is more natural to the construction of Gaussian beams. Second, the classification result of singular Lagrangian maps is more complete than that of singular exponential maps in the literature, though these two problems are equivalent in Riemannian manifold.
Finally, it is more natural to study caustics in the cotangent space.

\section{Statement of the main results}


\subsection{Sensitivity of recovering the scattering relation from the DDtN map}

It is known that the DDtN map $\Lambda_c$ determines the scattering relation $\mathfrak{S}_{c}$ uniquely \cite{PU06}.
We show that the following sensitivity result of recovering the scattering relation from the DDtN map holds. The proof is given
in Section \ref{section-proof-of-thm1}.

\begin{thm} \label{thm1}
Let $c$ and $\tilde{c}$ be two velocity fields in the class
$\mathfrak{A}(\epsilon_0, \Omega, M_0, T)$.  Then there exists a constant $\delta>0$ such that
$$
\mathfrak{S}_{\tilde{c}}= \mathfrak{S}_{c}
$$
if $\|\Lambda_{\tilde{c}}- \Lambda_{c}\|_{H_0^1[0, 3\epsilon_1/4]\times \Gamma \rightarrow L^2([0,T+\epsilon_1]\times \Gamma)} \leq \delta$.
\end{thm}

\begin{rmk}
The same result holds when the velocity fields are replaced by symmetric positive definite matrices.
\end{rmk}

\subsection{Linearization of the operator which maps velocity fields to Hamiltonian flows}

We begin with the following observation.
\begin{lem}
Let $c$ and $\tilde{c}$ be two velocity fields in the class $\mathfrak{A}(\epsilon_0, \Omega, M_0, T)$, then
$\mathfrak{S}_{c} = \mathfrak{S}_{\tilde{c}}$ if and only if
$\mathcal{H}_{c}^T|_{S^*_{-}\Gamma} = \mathcal{H}_{\tilde{c}}^T|_{S^*_{-}\Gamma}$.
\end{lem}

The above lemma shows the equivalence of the Hamiltonian flow and the scattering relation.
The next lemma shows that $\mathcal{H}_c^t$ satisfies an equivalent ordinary differential equation (ODE) system in $S^*\R^d$.

\begin{lem}
Let $(x_0, \xi_0)\in S^*\R^d =\{(x, \xi)\in \R^{2d}: |\xi|=1\}$, and let $(x(t), \xi(t))=\mathcal{H}_c^t(x_0, \xi_0)$, then $(x(t), \xi(t))$ satisfies the following ODE system
\begin{eqnarray}
  \dot{x} &=&\f{\xi}{\xi_1^2+...+\xi_d^2}, \label{h1}\\
  \dot{\xi} &=& b(x). \label{h2}
\end{eqnarray}
where $b(x)= -\f{1}{2}\nabla \ln{c^2}$. Conversely, if $(x(t), \xi(t))\in S^*\R^d$ satisfies the ODE system (\ref{h1})-(\ref{h2}), then
$(x(t), \xi(t))=\mathcal{H}_c^t(x_0, \xi_0)$.
\end{lem}

We next linearize the operator which maps each velocity field to its induced Hamiltonian flow restricted to the cosphere bundle.
Let $c$ be a fixed smooth background velocity field.
Denote the perturbed velocity field and Hamiltonian flow at time $T$ as $\tilde{c}^2= c^2 + \delta c^2$ and
$\mathcal{H}_{\tilde{c}}^{T}= \mathcal{H}_c^{T} + \delta \mathcal{H}_c^{T}$ respectively. Denote also that
$\delta b = -\f{1}{2}\nabla (\ln{\tilde{c}^2} -\ln{c^2})$
and
$$
A(x, \xi)=\left( \begin{array}{cc}
0 & \f{\p}{\p \xi}(\f{\xi}{\xi_1^2+...+\xi_d^2}) \\
\f{\p b}{\p x} &  0  \\
\end{array} \right).
$$
For each $(x_0, \xi_0)\in S_{-}^*\Gamma$, let $\Upsilon(t, x_0, \xi_0)$ be the solution of the following ODE system
$$
\dot{\Upsilon}(t)= -\Upsilon(t)A(\mathcal{H}_c^t(x_0, \xi_0)), \q \Upsilon(0)= Id.
$$

By the results in Appendix \ref{appendix-ode}, we have
$$\delta \mathcal{H}_c^T= \f{\delta \mathcal{H}_c^T}{\delta b}(\delta b) +r(\delta b),
$$
where
\be \label{tau1}
\f{\delta \mathcal{H}_c^T}{\delta b}(\delta b)(x_0, \xi_0)= \int_0^T\Upsilon^{-1}(T, x_0, \xi_0)\cdot \Upsilon(s, x_0, \xi_0)\left( \begin{array}{c}
0  \\
\delta b(x(s, x_0, \xi_0)) \\
\end{array} \right)\,ds
\ee
and $\|r(\delta b)\|_{L^{\infty}} \leq C \|\delta b\|_{C^1}^2$
for some constant $C>0$ depending only on $\|c\|_{C^3(\R^d)}$.

Formula (\ref{tau1}) motivates us to define the following geodesic X-ray transform operator
\be \label{tau}
\mathfrak{I}_c(f)(x_0, \xi_0)=  \int_0^T\Upsilon(s, x_0, \xi_0)
f(x(s, x_0, \xi_0))\,ds, \quad f\in \mathcal{E}'(\Omega, \R^{2d}).
\ee

Then $\f{\delta \mathcal{H}_c^T}{\delta b}(\delta b)(x_0, \xi_0)= \Upsilon^{-1}(T, x_0, \xi_0) \cdot \mathfrak{I}_c(f)(x_0, \xi_0)$ with
\be \label{f}
f=\left( \begin{array}{c} 0  \\
\f{1}{2}\nabla (\ln c^2- \ln \tilde{c}^2) \\
\end{array} \right).
\ee

We associate each $(x, \xi)\in S^*\Om$ a matrix $\Phi(x, \xi)$. Let
$(x_0, \xi_0)= \tau(x, \xi)=\mathcal{H}_c^{l_{-}(x, \xi)}(x, \xi)$.
We then define
$$
\Phi(x, \xi)= \Upsilon(-l_{-}(x, \xi),\tau(x, \xi)).
$$

By definition, it is clear that
$$
\Phi(\mathcal{H}_c^{s}(x_0, \xi_0))=\Upsilon(s, x_0, \xi_0)
$$
for all $s\in \R_{+}$ such that $\mathcal{H}_c^{s}(x_0, \xi_0)\in S^*\Om$.
Note that $f(x(s, x_0, \xi_0))=0$ if $s\geq l(x_0, \xi_0)$ for $f\in \mathcal{E}'(\Omega, \R^{2d})$.
We can rewrite the X-ray transform operator $\mathfrak{I}_c$ (\ref{tau}) in the following standard form
\bea
\mathfrak{I}_c f(x_0, \xi_0)
&=& \int_0^{l(x_0, \xi_0)}\Upsilon(s, x_0, \xi_0)
f(x(s, x_0, \xi_0))\,ds,  \nonumber\\
&=& \int_0^{l(x_0, \xi_0)}\Phi(\mathcal{H}_c^{s} (x_0, \xi_0))f(\pi(\mathcal{H}_c^{s} (x_0, \xi_0))) \,ds. \label{formula-I}
\eea

\begin{rmk}
Formula (\ref{formula-I}) is derived in the coordinate of $T^*\R^d$. Hence it may not be geometrically invariant.
\end{rmk}

\begin{lem} \label{lem-scattering-relation}
Assume that $\mathfrak{S}_{c} = \mathfrak{S}_{\tilde{c}}$, let $f$ be defined as in (\ref{f}), then
\[
 \|\mathfrak{I}_c f\|_{L^{\infty}} \lesssim \|f\|_{C^1(\Omega)}^2.
\]
\end{lem}

{\bf{Proof}}: First, using $\mathfrak{S}_{c} = \mathfrak{S}_{\tilde{c}}$ and Lemma 3.1, we have
\begin{equation} \label{eq-a1}
\delta \mathcal{H}_c^T|_{S^*_{-}\Gamma}=\mathcal{H}_{\tilde{c}}^T|_{S^*_{-}\Gamma}-\mathcal{H}_{c}^T|_{S^*_{-}\Gamma}= 0.
\end{equation}
Next, by the proceeding discussion, we have
\be \label{eq-a2}
\delta \mathcal{H}_c^T= \frac{\delta \mathcal{H}_c^T}{\delta b}(\delta b) +r(\delta b),
\ee
where the term $r$ satisfies the following inequality
\be \label{eq-a3}
\|r(\delta b)\|_{L^{\infty}} \leq C \|\delta b\|_{C^1}^2
\ee
for some constant $C>0$ depending only on $\|c\|_{C^3(\R^d)}$.
Combining (\ref{eq-a1})-(\ref{eq-a3}), we see that
$$
\|\frac{\delta \mathcal{H}_c^T}{\delta b}(\delta b)\|_{L^{\infty}} = \|r(\delta b)\|_{L^{\infty}} \lesssim  \|\delta b\|_{C^1}^2 =\|f\|_{C^1}^2.
$$
Finally, using the equality
$$\mathfrak{I}_c(f)(x_0, \xi_0)= \Upsilon(T, x_0, \xi_0)\cdot \frac{\delta \mathcal{H}_c^T}{\delta b}(\delta b)(x_0, \xi_0)$$
and the fact that the matrix-valued function $\Upsilon(T, \cdot, \cdot)$ is a smooth function determined by the background velocity field $c$, we get the desired conclusion immediately.

\subsection{Fold-regular points and local stability for geodesic X-ray transform} \label{subsec-fold-regular}
We consider the stability estimate of the operator $\mathfrak{I}_c $. For simplicity, we drop the subscript $c$.
Define $\beta: T^*\R^d \backslash \{(x, 0): x\in \R^d\} \to S^*\R^d$ by
\be  \label{beta}
 \beta (x, \xi )= \left(x, \f{\xi}{|\xi|}\right).
\ee

Let $\pi: T^*\R^d \to \R^d$ be the natural projection onto the base space. We define $\phi: T^*\R^d \to \R^d$ by
$$
\phi(x, \xi)= \pi \circ \mathcal{H}^{t=1}(x, \xi), \quad (x, \xi)\in T^*\R^d.
$$
We remark that $\phi$ defined above is equivalent to the exponential map in Riemannian manifold.

The following result about the normal operator $\mathfrak{N}= \mathfrak{I}^{\dag}\mathfrak{I}$ is well-known.

\begin{lem} \label{lem-x-ray3}
The normal operator $\mathfrak{N}: L^2(\Omega, \R^{2d}) \to L^2(\Omega, \R^{2d})$ is bounded and has the following representation
\be \label{normal-operator}
\mathfrak{N}f(x)= \int_{T^*_{x}\Omega}W(x, \xi)f(\phi(x, \xi))\,d\sigma_x(\xi), \quad f \in L^2(\Omega, \R^{2d})
\ee
where $d\sigma_x$ denotes
the measure in the space $T_x^*\R^d$ induced by the velocity field $c$, i.e. $d\sigma_x(\xi)= c(x)^{d}d\xi$, and $W$ is defined as
\be \label{formula-w}
W(x, \xi)= \f{1}{|\xi|^{d-1}}\{ \Phi^\dag \circ \beta(x, \xi) \cdot \Phi\circ \beta\circ \mathcal{H}(x, \xi)
+ \Phi^\dag\circ \beta(x, -\xi) \cdot \Phi\circ \beta\circ \mathcal{H}^{-1}(x, -\xi) \}.
\ee
\end{lem}
{\bf{Proof}}. See \cite{SUAPDE} or \cite{SU04Duke}.

\medskip

We see from (\ref{normal-operator}) that the local property of the normal operator $\mathfrak{N}$
restricted to a small neighborhood of $x \in \Omega$
is determined by the lagrangian map $\phi(x, \cdot): T^*_{x}\R^d \to \R^d$. When the map is a diffeomorphism, it is known that the operator
$\mathfrak{N}$ near $x$ is a pseudo-differential operator ($\Psi$DO). However, in general case, the map may not be a diffeomorphism and may have singular points which are called caustic vectors. The value of the map at caustic vectors are called caustics.
When caustics occur, the Schwartz kernel of the operator
$\mathfrak{N}$ has two singularities, one is from the diagonal which contributes to a $\Psi$DO $\mathfrak{N}_1$, and the other is from the caustics which contributes to a singular integral operator $\mathfrak{N}_2$. The property of $\mathfrak{N}_2$ depends on the type of caustics. The case for fold caustics is investigated in \cite{SUAPDE}, where it is shown that fold caustics contribute a Fourier Integral Operator (FIO) to $\mathfrak{N}_2$. Little is known for caustics of other type.
Here we recall the following definition of fold caustics.
\begin{definition}
Let $f: \R^n \to \R^n$ be a germ of $C^{\infty}$ map at $x_0$, then $x_0$ is said to be a fold vector and $f(x_0)$ a fold caustic
if the following two conditions are satisfied:
\begin{enumerate}
\item
the rank of $df$ at $x_0$ equals to $n-1$ and $\det df$ vanishes of order 1 at $x_0$;

\item
the kernel of the matrix $df(x_0)$ is transversal to the manifold $\{x: \det df(x)=0\}$ at $x_0$.
\end{enumerate}
\end{definition}

We now introduce the following concept of ``operator germ'' to characterize the contribution of an infinitesimal neighborhood
of a caustic or a regular point to the normal operator $\mathfrak{N}$.
\begin{definition}
For each $\xi \in T^*_{x}\R^d \backslash 0$, the operator germ $\mathfrak{N}_{\xi}$ is defined to be the equivalent class of operators in the following form
\be \label{operator-germ}
\mathfrak{N}_{\xi} f(y)= \int_{T^*_{y}\Omega}W(y, \eta)f(\phi(y, \eta)) \chi(y, \eta)\,d\sigma_y(\eta).
\ee
where $\chi$ is a smooth function supported in a small neighborhood of $(x, \xi)$ in $\R^{2d}$.
Two operators with $\chi_1$ and $\chi_2$
are said to be equivalent if there exists a neighborhood $B(x, \xi)$ of $(x, \xi)$
such that $\chi_1=\chi_2\cdot \chi_3$ for some $\chi_3 \in C^{\infty}_0(B(x, \xi))$ with $\chi_3(x, \xi) \neq 0$.

The operator germ $\mathfrak{N}_{\xi}$ is said to has certain property if there exists a neighborhood $B(x, \xi)$ of $(x, \xi)$ in $T^*\R^d$ such that the property holds for all operators of the form (\ref{operator-germ}) with $\chi \in C^{\infty}_0(B(x, \xi))$.
\end{definition}

Properties of the above defined operator germ will be given in Section \ref{subsec-operator-germ}.

We note from the preceding discussion that it is complicated to
analyze the full operator $\mathfrak{N}$ which contains
information from all geodesics. However, for a given interior point $x$, to recover $f$ or the singularity of $f$
at $x$ from its geodesic transform, we need only to
select a set of geodesics whose conormal bundle can cover the cotangent space $T_x^*\R^d$. Caustics may be allowed along these geodesics as long as  they are of the simplest type, i.e fold type so that we can analyze their contributions. This idea can be carried out by introducing a cut-off function for the set of geodesics as we do now. We remark that this idea is motivated by the work \cite{SU08AJM}.
For any $\alpha\in C^{\infty}_0(S^*_{-}\Gamma)$, we define
\begin{equation} \label{formula-I-alpha}
\mathfrak{I}_{\alpha}f(x_0, \xi_0)
= \alpha(x_0, \xi_0) \int_0^{l(x_0, \xi_0)}\Phi(\mathcal{H}_c^{s} (x_0, \xi_0))f(\pi(\mathcal{H}_c^{s} (x_0, \xi_0))) \,ds
\end{equation}
where $(x_0, \xi_0) \in S^*_{-}\Gamma$.
Let $\alpha^{\sharp}$ be the unique lift of $\alpha$ to $S^*\Omega$ which is constant along bicharacteristic curves, i.e.
$\alpha^{\sharp}(x, \xi) = \alpha\circ \tau(x, \xi)$ for $(x, \xi) \in S^*\Omega$.
Then $\alpha^{\sharp}$ is smooth in $S^*\Omega$ and we have
\be \label{formula-I-alpha-1}
\mathfrak{I}_{\alpha}f(x_0, \xi_0) = \int_0^{l(x_0, \xi_0)} (\alpha^{\sharp}\cdot \Phi)(\mathcal{H}_c^{s} (x_0, \xi_0))f(\pi(\mathcal{H}_c^{s} (x_0, \xi_0))) \,ds
\ee

With the original weight $\Phi$ being replaced by the new one $\alpha^{\sharp}\cdot \Phi$, we similarly can define
$\mathfrak{N}_{\alpha}$. In fact, it is easy to check that $\mathfrak{N}_{\alpha}$ is defined as in (\ref{formula-w}) with $W$ being replaced by
\beas
W_{\alpha}(x, \xi) &=&  \f{1}{|\xi|^{d-1}}
  |\alpha\circ \tau\circ \beta(x, \xi)|^2\Phi^\dag \circ \beta(x, \xi) \cdot \Phi\circ \beta\circ \mathcal{H}(x, \xi)\\
 && + \f{1}{|\xi|^{d-1}}|\alpha\circ \tau\circ \beta(x, -\xi)|^2\Phi^\dag \circ \beta(x, -\xi) \cdot \Phi\circ \beta\circ \mathcal{H}^{-1}(x, -\xi).
\eeas

It can be shown that with properly chosen $\alpha$, the analysis of the operator $\mathfrak{N}_{\alpha}$ becomes possible and we can recover singularities of $f$ from $\mathfrak{N}_{\alpha}f$.

We now give two definitions whose discussions are postponed to Section 6.
\begin{definition}
A fold vector $\xi\in T^*_{x}\R^d$ is called fold-regular if there exists a neighborhood $U(x)$ of $x$ such that
the operator germ $\mathfrak{N}_{\xi}$
is compact from
$L^2(\Omega_{\epsilon_0}, \R^{2d})$ to $H^1(U(x), \R^{2d})$ (or from
$H^s(\Omega_{\epsilon_0}, \R^{2d})$ to $H^{s+1}(U(x), \R^{2d})$ for all $s\in \R$).
\end{definition}

\begin{definition} \label{def-fold-regular-point}
A point $x$ is called fold-regular if there exists a compact subset $\mathcal{Z}(x) \subset S^*_{x}\R^d$ such that the following two conditions are satisfied:
\begin{enumerate}
\item
For each $\xi \in \mathcal{Z}(x)$, there exist either no singular vectors or singular vectors of fold-regular type along the ray $\{t\xi:\, t\in \R\}$ for the map
$\phi(x, \cdot)$;

\item
$\forall \, \xi \in S^*_{x}\R^d, \, \exists \, \theta \in \mathcal{Z}(x) , \mbox{ such that }\, \theta \perp \xi$.
\end{enumerate}
\end{definition}

We remark that $\mathcal{Z}(x)$ parameterizes a subset of geodesics that pass through $x$ and along which there exist either no caustics or caustics of fold-regular type. A set of geodesics satisfying condition 2 in the above definition is called complete.

We now present the main result on the local stability estimate for the geodesic X-ray transform operator. The proof is given in Section 6.
\begin{thm} \label{thm2}
Let $x_*$ be a fold-regular point, then there exist a cut-off function $\alpha\in C^{\infty}_0(S^*_{-}\Gamma)$, a neighborhood $U(x_*)$ of $x_*$,
a compact operator $\mathfrak{N}_{2, \alpha}$ from
$L^2(\Omega_{\epsilon_0}, \R^{2d})$ to $H^1(U(x_*), \R^{2d})$ and
a smoothing operator $\mathfrak{R}$ from $\mathcal{E}'(\Omega, \R^{2d})$ into $C^{\infty}(\overline{U(x_*)}, \R^{2d})$,
such that for any $U_0(x_*) \Subset U(x_*)$ the following holds
\begin{equation}{\label{11}}
\|f\|_{H^s(U_0(x_*), \R^{2d})} \lesssim   \|\mathfrak{N}_{\alpha}f \|_{H^{s+1}(U(x_*), \R^{2d})}+ \|\mathfrak{N}_{2, \alpha}f \|_{H^{s+1}(U(x_*), \R^{2d})}+
\|\mathfrak{R} f\|_{H^{s}(U(x_*), \R^{2d})}
\end{equation}
for all $f\in \mathcal{D}'(\Omega_{\epsilon_0}, \R^{2d})$ and $s\in \R$.
\end{thm}

\subsection{Sensitivity of recovering the velocity field from the DDtN map} \label{subsec-thm3}
\begin{definition}
An admissible velocity field $c$ is called fold-regular if all points in $\Omega$ are fold-regular with respect to the Hamiltonian flow $\mathcal{H}_c^t$.
\end{definition}

We have established the following main result on the sensitivity of recovering velocity field from DDtN map.
For simplicity we only consider the case $d=3$, similar results also hold for $d>3$.
The proof is given in Section 7.

\begin{thm} \label{thm3}
Let $c$ and $\tilde{c}$ be two velocity fields in the class
$\mathfrak{A}(\epsilon_0, \Omega, M_0, T)$.  Assume that the velocity field $c$ is smooth and is fold-regular.
Then there exist a finite dimensional subspace $\mathfrak{L}\subset L^2(\Omega_{\epsilon_0}, \R^{3})$,
and a constant $\delta> 0$ such that for all
$\tilde{c}$ sufficiently close to $c$ in $H^{\f{17}{2}}(\Omega)$ and satisfying $\nabla (\ln c^2- \ln \tilde{c}^2) \perp \mathfrak{L}$, $\|\Lambda_{\tilde{c}}- \Lambda_{c}\|_{H_0^1[0, 3\epsilon_1/4]\times \Gamma \rightarrow L^2([0,T+\epsilon_1]\times \Gamma)} \leq \delta$
implies that $c=\tilde{c}$.
\end{thm}

\section{Gaussian beam solutions to the wave equation}
Let $c$ be a velocity field in the class $\mathfrak{A}(\epsilon_0, \Omega, M_0, T)$.
We construct Gaussian beam solutions to the wave equation system (\ref{wave1})-(\ref{wave3}) in this section.

We first construct a Gaussian beam in $\R^d$. Following \cite{QY1}, we define
$G(x, \xi)= |\xi|=c(x)\sqrt{\xi_1^2+...+\xi_d^2}$. For a given $(x_0, \xi_0)\in S^*_{-}\Gamma$, let $(x(t), \xi(t), M(t), a(t))$ be the solution to the following ODE system:
\bea
  \dot{x} &=& G_p, \quad \q \q x(t_0) = x_0, \nonumber\\
  \dot{\xi} &=& -G_x, \quad \q \q \xi(t_0)= {\xi}_0, \nonumber\\
  \dot{M} &=& -G_{x \xi}^{\dag} M - M G_{\xi x} - M G_{\xi \xi}M - G_{xx}, \quad M(t_0) = \sqrt{-1}\cdot Id, \label{ode3}\\
  \dot{a} &=& -\frac{a}{2G}(c^2\mbox{trace}(M)-G_x^\dag G_{\xi}-G_{\xi}^{\dag}MG_{\xi}),\quad a(t_0) =\lambda^{\f{d}{4}}. \label{transport}
\eea

The corresponding Gaussian beam with frequency $\lambda$ ($\lambda \gg 1$) is given as follows
$$
g(t,x, \lambda)= a(t)e^{i \lambda \tau(t,x)}
$$
where $\tau(t,x)= \xi(t)\cdot(x-x(t))+ \f{1}{2}(x-x(t))^{\dag}M(t)(x-x(t))$.

\medskip

Now, let the beam $g$ impinge on the surface $\Gamma$ transversely, we want to construct the reflected beam $g^-$. Let
the ray $x(t)$ hits $\Gamma$ at the point $x(t_1)=x_1$. Write $\xi(t_1)=\xi_1$. We parameterize $\Gamma$ in a neighborhood of $x_1$, say $V(x_1)$, by a smooth diffeomorphism $F: U(x_1) \to V(x_1)$, where $U(x_1)$ is a neighborhood of the origin in $\R^{d-1}$. We require that $F(0)=x_1$.
With the coordinate $x=F(y)$, we can rewrite functions when their spatial variables are restricted to the boundary $\Gamma$. For example,
we rewrite
$$
g(t,x)= g(t, F(y))= \hat{g}(t, y), \quad \tau(t,x)= \tau(t,F(y))= \hat{\tau}(t, y), \quad \mbox{for } x \in V(x_1).
$$
We remark that here and throughout this section and the next, unless specified otherwise, we use the notation $\hat{f}$ to denote the function $f$ restricted to the space-time boundary $\Gamma \times [0, T]$ under the coordinate $x=F(y)$ for the spatial variables.

By direct calculation, we have
\beas
\hat{\tau}(t_1, y)
&=& \hat{\tau}(t_1, 0)+ (\f{\p \hat{\tau}}{\p (t, y)}(t_1, 0))^{\dag}(t-t_1, y)+
(t-t_1, y)\f{\p^2 \hat{\tau}}{\p (t, y)^2}(t_1, 0)(t-t_1, y)^{\dag} \\
&& + O(|(t-t_1, y)|^3) \\
&=&\langle (-1, \f{\p F}{\p y}(0)^{\dag}\xi_1), (t-t_1, y)\rangle + (t-t_1, y)\hat{M}(t_1)(t-t_1, y)^{\dag}+ O(|(t-t_1, y)|^3),
\eeas
where the matrix $\hat{M}(t_1)$ is defined as the Hessian of the phase $\hat{\tau}$ at $(t_1, 0)$, i.e. $\hat{M}(t_1)=\f{\p^2 \hat{\tau}}{\p (t,y)^2}(t_1, 0)$. It is clear that $\hat{M}(t_1)$ is determined by $M(t_1)$ and the coordinate function $x=F(y)$. Using the assumption that the ray $x(t)$ intersects $\Gamma$ transversely, we can conclude that $\Im{\hat{M}(t_1)}>0$. Moreover, if condition (\ref{condition-ad2}) is satisfied, we have
\be  \label{imginary-part}
\Im{\hat{M}(t_1)}>C
\ee
for some $C>0$ depending only on $\epsilon_0$ and $M_0$.

\medskip

We proceed to construct the reflected beam $g^-$. Write
$$
g^{-}(t,x, \lambda)= a^{-}(t)e^{i \lambda \tau^{-}(t,x)}
$$
with
$$
\tau^{-}(t,x)= \xi^{-}(t)\cdot(x-x^{-}(t))+ \f{1}{2}(x-x^{-}(t))^{\dag}M^{-}(t)(x-x^{-}(t)).
$$
We need to find $(x^{-}(t_1), \xi^{-}(t_1), a^{-}(t_1), M^{-}(t_1))$ such that the $g^{-}+g \thickapprox 0$ on the boundary. Following \cite{AAB}, we impose the following condition
\be \label{conditon-tau}
\p^{\alpha}_{t, y}\hat{\tau}(t_1, 0)= \p^{\alpha}_{t, y}\hat{\tau}^{-}(t_1, 0), \quad \mbox{for all } \, |\alpha|\leq 2.
\ee
As a result, we obtain $\xi^{-}(t_1)= \xi^{-}_1= \xi_1- 2 \langle \xi_1, \nu(x_1) \rangle \nu(x_1)$ and
$\hat{M}^-(t_1)= \hat{M}(t_1)$. Here we note that $\hat{M}^{-}(t_1)$ is the Hessian of the phase $\hat{\tau}^{-}$ for the reflected beam at the point $(t_1, 0)$.
Consequently, $M^{-}(t_1)$ is also determined.
Finally, set $x^{-}(t_1)=x_1$ and $a^{-}(t_1)= -a(t_1)$. Then all of the four components of $(x^{-}(t_1), \xi^{-}(t_1), a^{-}(t_1), M^{-}(t_1))$ are constructed.
 Afterward, we solve an ODE system to get  $(x^{-}(t), \xi^{-}(t), M^{-}(t), a^{-}(t))$ as we did for the beam $g$. This completes the construction for the reflected beam $g^-$.
\medskip

We now present some properties about the constructed beam. The following lemma is crucial in the subsequent estimates. We refer to \cite{R83} for the proof.
\begin{lem} \label{lem-phase-matrix}
Both the matrices $M(t)$ and $M^{-}(t)$ are uniformly bounded for $t\in [0, T+ \epsilon_1]$. Moreover, there exists $C>0$, depending on $M_0$ and $\epsilon_0$, such that $\Im{M}(t)>C$ and $\Im{M}^{-}(t)>C$ for all $t\in [0, T+ \epsilon_1]$.
\end{lem}

We next introduce two auxiliary beams $\hat{g}_*(t, y, \lambda) = a(t_1) e^{i \lambda \hat{\tau}_*}$ and $
\hat{g}^-_*(t, y, \lambda) =a^{-}(t_1) e^{i \lambda \hat{\tau}_*^{-}}$,
where
\beas
\hat{\tau}_* &=&   \bigl\langle (-1, \f{\p F}{\p y}(0)^{\dag}\xi_1), (t-t_1, y)\bigr\rangle + (t-t_1, y)\hat{M}(t_1)(t-t_1, y)^{\dag}, \\
\hat{\tau}^-_* &=&  \bigl \langle (-1, \f{\p F}{\p y}(0)^{\dag}\xi^-_1), (t-t_1, y)\bigr\rangle + (t-t_1, y)\hat{M}^{-}(t_1)(t-t_1, y)^{\dag}.
\eeas

Compared to the Gaussian beams $\hat{g}$ and $\hat{g}^{-}$, the axillary beams $\hat{g}_*$ and $\hat{g}^-_*$ have frozen amplitude at $t=t_1$ and phase function with only quadratic terms.
From the proceeding construction of the reflected beam, it is clear that $\hat{\tau}_*= \hat{\tau}^-_*$ and $\hat{g}_*= -\hat{g}^-_*$. Moreover, as is shown in the next lemma, Lemma \ref{lem-g-ghat1}, the two axillary beams $\hat{g}_*$ and $\hat{g}^-_*$ are good approximations to the incident beam $\hat{g}$ and reflected beam $\hat{g}^{-}$ on the space-time boundary respectively. This property is used for estimating the interactions of Gaussian beams on the space-time boundary (see Step 3 in the proof of Theorem \ref{thm1}).

\begin{lem} \label{lem-g-ghat1}
\bea
\hat{g}(t, y, \lambda) &=& \hat{g}_*(t, y,  \lambda) + O(\sqrt{\lambda}) \quad \mbox{in }\,\, H^1((3\epsilon_1/4, t_1+ \epsilon_1/2)\times U(x_1)), \label{g-approximation1-1} \\
\hat{g}^-(t, y,  \lambda)
&=&  \hat{g}_*^{-}(t, y,  \lambda) + O(\sqrt{\lambda}) \quad \mbox{in }\,\, H^1((3\epsilon_1/4, t_1+ \epsilon_1/2)\times U(x_1)).
\label{g-approximation2-2}
\eea
\end{lem}

{\bf{Proof:}} See Appendix \ref{appendix-proof}.

\medskip

Note that $\|\hat{g}(t, y, \lambda)\|_{L^2((t_1-\epsilon_1/2, t_1+\epsilon_1/2)\times U(x_1))}  \thickapprox  1$.
As a direct consequence of
Lemma \ref{lem-g-ghat1}, we obtain the following norm estimate for the beam $g$ restricted to the boundary $\Gamma$.

\begin{lem} \label{g-norm}
\be
\|g(\cdot, \cdot, \lambda)\|_{L^2((t_1-\epsilon_1/2, t_1+\epsilon_1/2)\times V(x_1))}  \thickapprox  1.
\ee
\end{lem}

We now present an $H^1$-norm estimate for $g^- + g$ and an approximation for
the Neumann data $\f{\p g}{\p \nu}^{-}+ \f{\p g}{\p \nu}$ on the boundary.
\begin{lem}  \label{lem-g-sum}
\bea
 g^{-}(t, x, \lambda)+  g(t, x, \lambda) &=&
O(\sqrt{\lambda}) \quad \mbox{in \,\,} H^1((3\epsilon_1/4, t_1+ \epsilon_1/2)\times V(x_1)); \label{gb-boundary1}\\
 \f{\p g}{\p \nu}^{-}+ \f{\p g}{\p \nu} &=& 2 i \lambda g \cdot \langle \xi_1, \nu(x_1)\rangle +
O(\sqrt{\lambda}) \,\, \mbox{in \,} L^2((3\epsilon_1/4, t_1+ \epsilon_1/2)\times V(x_1))  \label{gb-boundary2}
\eea
\end{lem}

{\bf{Proof:}} See Appendix \ref{appendix-proof}.

\medskip

Now, we are ready to construct Gaussian beam solutions to the initial boundary value problem of the wave system (\ref{wave1})-(\ref{wave3}).
We first choose $\chi_{\epsilon_1}(t) \in C_0^{\infty}(\R)$ such that $\chi_{\epsilon_1}(t)=1$ for $t\in (\epsilon_1/4, \epsilon_1/2)$ and
$\chi_{\epsilon_1}(t)=0$ for $t\in (-\infty, 0) \bigcup (3\epsilon_1 /4 , \infty)$.
Let $(x_0, \xi_0)\in S^*_{-}\Gamma$ and
$(x_0^*, \xi_0^*)= \mathcal{H}^{-\f{\epsilon_1}{4}}(x_0, \xi_0)=(x_0 - \f{\epsilon_1 \cdot \xi_1}{4}, \xi_0)$.  Let $g$ be the Gaussian beam constructed with the initial data $x(0)= x_0^*, \xi(0)= \xi_0^*, M(0)= i\cdot Id$ and $a(0)= \lambda^{\f{d}{4}}$.
The beam $g$ is reflected by $\Gamma$ at $(x_1, \xi_1)= \mathfrak{S}_{c}(x_0, \xi_0)= \mathcal{H}_{c}^{l(x_0, \xi_0)}(x_0, \xi_0)$ at
$t_1= l(x_0, \xi_0)+ \f{\epsilon_1}{4}$. We construct
the reflected beam $g^{-}$ by the preceding procedure.
Let $u$ be the exact solution to the wave system (\ref{wave1})-(\ref{wave3}) with
$$
f(t, x, \lambda)= g(t, x, \lambda)\cdot \chi_{\epsilon_1}(t).
$$

Then $u=g+g^- + R$, where the remaining term $R$ satisfies the following equation system
\beas
\mathcal{P} R&=&- \mathcal{P}(g+g^-),   \quad (t, x)\in \Om \times(0, t_1+\epsilon_1/2), \\
R(0,x, \lambda)&=&- (g+g^-)(0, x, \lambda),                  \quad  x\in \Omega , \\
R_t(0,x, \lambda)&=&-(g_t+ g^-_t)(0, x, \lambda),                  \quad  x\in \Omega , \\
R(t,x, \lambda)&=&- g(t, x, \lambda)(1-\chi_{\epsilon_1}(t))- g^-(t, x, \lambda),           \quad   (t, x)\in (0,t_1+\epsilon_1/2) \times \Gamma .
\eeas
Here $\mathcal{P}$ stands for the wave operator $\frac{1}{c^2(x)}\p_{tt}-\Delta $.

\begin{lem} \label{reminder-estimate}
$$
\|\f{\p R}{\p \nu}\|_{L^2([0, t_1+ \epsilon_1/2]\times \Gamma)}
\leq C \sqrt{\lambda}
$$
for some constant $C>0$ depending on $\epsilon_0$ and $M_0$.
\end{lem}

{\bf{Proof}}.  We apply Theorem 4.1 in \cite{LLT86} to derive the estimate. Note that the compatibility condition is satisfied on the boundary at time $t=0$. It remains to show that the following four estimates hold:
\bea
\|\mathcal{P}(g+g^-)\|_{C([0,t_1+ \epsilon_1/2]; L^2(\Omega) )} & \lesssim & \sqrt{\lambda}, \label{estimate1}\\
\|(g+g^-)(0, \cdot, \lambda)\|_{H^1(\Omega)} & \lesssim & \sqrt{\lambda}, \label{estimate2}\\
\|(g_t+g^-_{t})(0, \cdot, \lambda)\|_{L^2(\Omega)} & \lesssim & \sqrt{\lambda}, \label{estimate3}\\
\|g(t, x, \lambda)(1-\chi_{\epsilon_1}(t))- g^-(t, x, \lambda)\|_{H^1([0, t_1+ \epsilon_1/2]\times \Gamma)} & \lesssim & \sqrt{\lambda}.\label{estimate4}
\eea
First, (\ref{estimate1}) follows from the standard estimate for Gaussian beams, see for example \cite{BQ1}.
We next show (\ref{estimate2}).  By Lemma \ref{lem-phase-matrix}, there exists a constant $C>0$ depending on
$M_0$ and $\epsilon_0$ such that
the following two inequalities hold
\beas
|g(t, x, \lambda)| & \lesssim & \lambda^{\f{d}{4}}\cdot e^{-C\lambda\cdot |x-x^{-}(t)|^2},\\
|g^{-}(t, x, \lambda)| & \lesssim & \lambda^{\f{d}{4}}\cdot e^{-C\lambda\cdot |x-x^{-}(t)|^2}.
\eeas
Thus the beam $g$ and $g^{-}$ are exponentially decaying away from the ray $x(t)$ and $x^{-}(t)$ respectively.
Using this property, it is straightforward to show that
$\|g(0, \cdot, \lambda)\|_{H^1(\Omega)} \lesssim  1$ and $\|g^-(0, \cdot, \lambda)\|_{H^1(\Omega)} \lesssim  1$,
whence (\ref{estimate2}) and (\ref{estimate3}) follows.

Now, we show (\ref{estimate4}). We divide the domain $(0,t_1+\epsilon_1/2) \times \Gamma$ into three parts:
$$ \Sigma_1= (0, \epsilon_1/2) \times \Gamma,
\quad \Sigma_2= (\epsilon_1/2, t_1- \epsilon_1/2) \times \Gamma,
\quad
\Sigma_3= (t_1- \epsilon_1/2, t_1+ \epsilon_1/2) \times \Gamma.
$$
We show that inequality (\ref{estimate4}) holds on each part.

For $(t, x) \in \Sigma_1$, we have $1-\chi_{\epsilon_1}(t)=0$. Consequently,
$$
g(t, x)(1-\chi_{\epsilon_1}(t))- g^-(t, x, \lambda) = g^-(t, x, \lambda).
$$
By the exponential decaying property of $g^{-}$, we obtain that
$$
\|g(t, x)(1-\chi_{\epsilon_1}(t))- g^-(t, x, \lambda)\|_{H^1({\Sigma_1})} \lesssim \sqrt{1}.
$$
For $(t, x) \in \Sigma_2$,  by the exponential decaying property for both $g$ and $g^{-}$ again, we obtain
$$
\|g(t, x, \lambda)(1-\chi_{\epsilon_1}(t))- g^-(t, x, \lambda)\|_{H^1({\Sigma_2})} \lesssim  1.
$$
Finally, for $(t, x) \in \Sigma_3$, note that $t_1-\f{\epsilon_1}{2} = l(x_0, \xi_0)+ \f{\epsilon_1}{4}- \f{\epsilon_1}{2}\geq \f{3\epsilon_1}{4}$.
We can apply Lemma \ref{lem-g-sum} to the part $x\in V(x_1)$ and
the exponential decaying property for both $g$ and $g^{-}$ to the remaining part to conclude that
$$
\|g(t, x, \lambda)(1-\chi_{\epsilon_1}(t))- g^-(t, x, \lambda)\|_{H^1({\Sigma_3})} \lesssim  \sqrt{\lambda}
$$
This completes the proof of (\ref{estimate4}) and hence the lemma.


\section{Proof of Theorem \ref{thm1}} \label{section-proof-of-thm1}


\textbf{Proof of Theorem \ref{thm1}}. For any $(x_0, \xi_0)\in S^*_{-}\Gamma$, let
$(x_1, \xi_1)= \mathfrak{S}_{c}(x_0, \xi_0)= \mathcal{H}_{c}^{l(x_0, \xi_0)}(x_0, \xi_0)$ and
$(\tilde{x}_1, \tilde{\xi}_1)= \mathfrak{S}_{\tilde{c}}(x_0, \xi_0)= \mathcal{H}_{\tilde{c}}^{\tilde{l}(x_0, \xi_0)}(x_0, \xi_0)$.
We need to show that $(l(x_0, \xi_0), x_1, \xi_1)= (\tilde{l}(x_0, \xi_0), \tilde{x}_1, \tilde{\xi}_1)$ if
$\|\Lambda_{\tilde{c}}- \Lambda_{c}\|$ is sufficiently small.
We do this in the following steps.

\medskip
Step 1. Let $t_1= l(x_0, \xi_0)+ \f{\epsilon_1}{4}$ and
$\tilde{t}_1= \tilde{l}(x_0, \xi_0)+ \f{\epsilon_1}{4}$.
Without loss of generality, we may assume that $t_1\leq \tilde{t}_1$. Let $V(x_1)$ be a neighborhood of $x_1$ in $\Gamma$ which is parameterized by
a smooth function $F: U(x_1) \to V(x_1)$ as before. Note that if $\tilde{x}_1$ does not belong to $V(x_1)$, then we can show that $\|\Lambda_{\tilde{c}}- \Lambda_{c}\|$ is bounded from below by some positive constant and hence the result of Theorem \ref{thm1} is obvious. Therefore we may assume that $\tilde{x}_1 \in V(x_1)$. Let $\tilde{x}_1= F(\delta y)$.
We construct the initial beam $g$, the reflected beam $g^{-}$, the boundary Dirichlet data $f$, the solution $u$ to the wave equation with velocity field $c$ and remanning term $R$ as in the previous section. We similarly construct $\tilde{g}$, $\tilde{g}^-$, $\tilde{u}$ and $\tilde{R}$ to the system with velocity field $\tilde{c}$ and with boundary Dirichlet data $\tilde{f}=f$.

\medskip
Step 2. Denote by $I(t_1, \epsilon_1/2)$ the interval $(t_1-\epsilon_1/2, t_1+\epsilon_1/2 )$. Since $t_1\leq \tilde{t}_1$ and
$l(x_0, \xi_0)\geq \epsilon_1$, we have $I(t_1, \epsilon_1/2) \subset (3\epsilon_1/4, t_1+ \epsilon_1/2)$ and
$I(t_1, \epsilon_1/2) \subset (3\epsilon_1/4, \tilde{t}_1+ \epsilon_1/2)$. Then we can apply
(\ref{gb-boundary2}) and Lemma \ref{reminder-estimate} to obtain
\beas
(\Lambda_{\tilde{c}}- \Lambda_{c})f &=& \f{\p u}{\p \nu} - \f{\p \tilde{u}}{\p \nu} \\
&=& \f{\p (g+g^-)}{\p \nu} -\f{\p (\tilde{g}+\tilde{g}^-)}{\p \nu} + \f{\p R}{\p \nu} - \f{\p \tilde{R}}{\p \nu}\\
&=& 2i\lambda\cdot \bigl\{\langle \xi_1, \nu(x_1)\rangle \cdot g - \langle \tilde{\xi}_1, \nu(\tilde{x}_1)\rangle \cdot \tilde{g}\bigr\}
+ O(\sqrt{\lambda})
\eeas
in $L^2(I(t_1, \epsilon_1/2)\times V(x_1))$.

It follows that
\begin{align*}
\bigl\langle (\Lambda_{\tilde{c}}- \Lambda_{c})f, g \bigr\rangle_{L^2(I(t_1,\epsilon_1/2)\times V(x_1))} = 2i\lambda &\cdot \biggl[ \bigl\langle \xi_1, \nu(x_1)\bigr\rangle \cdot \bigl\langle g, g \bigr\rangle_{L^2(I(t_1,\epsilon_1/2)\times V(x_1))} \\ &-
\bigl\langle \tilde{\xi}_1, \nu(\tilde{x}_1)\bigr\rangle \cdot \bigl\langle \tilde{g}, g \bigr\rangle_{L^2(I(t_1,\epsilon_1/2)\times V(x_1))}\biggr]+ O(\sqrt{\lambda}).
\end{align*}
Note that
\beas
|\langle (\Lambda_{\tilde{c}}- \Lambda_{c})f, g \rangle_{L^2(I(t_1,\epsilon_1/2)\times V(x_1))}|
&\leq & \|(\Lambda_{\tilde{c}}- \Lambda_{c})f\|_{L^2(I(t_1,\epsilon_1/2)\times V(x_1))}\cdot \|g\|_{L^2(I(t_1,\epsilon_1/2)\times V(x_1))} \\
&\leq & \|(\Lambda_{\tilde{c}}- \Lambda_{c})f\|_{L^2((0,T+\epsilon_1)\times \Gamma)}\cdot \|g\|_{L^2(I(t_1,\epsilon_1/2)\times V(x_1))} \\
& \leq  & \|\Lambda_{\tilde{c}}- \Lambda_{c}\|_{H_0^1([0, 3\epsilon_1/4]\times \Gamma)\rightarrow L^2([0,T+\epsilon_1]\times \Gamma)}
\cdot \|f\|_{H^1_0([0,  3\epsilon_1/4]\times \Gamma)}\\
&& \cdot \|g\|_{L^2(I(t_1,\epsilon_1/2)\times V(x_1))} \\
& \lesssim  &
\lambda \cdot \|\Lambda_{\tilde{c}}- \Lambda_{c}\|_{H_0^1([0, 3\epsilon_1/4]\times \Gamma)\rightarrow L^2([0,T+\epsilon_1]\times \Gamma)}.
\eeas
Thus the following inequality holds
\bea
&& |\langle \xi_1, \nu(x_1)\rangle \cdot \langle g, g \rangle_{L^2(I(t_1,\epsilon_1/2)\times V(x_1))}|-
|\langle \tilde{\xi}_1, \nu(\tilde{x}_1)\rangle \cdot \langle \tilde{g}, g \rangle_{L^2(I(t_1,\epsilon_1/2)\times V(x_1))}| \nonumber\\
&& \leq \f{C}{\sqrt{\lambda}}+
C \cdot \|\Lambda_{\tilde{c}}- \Lambda_{c}\|_{H_0^1([0, 3\epsilon_1/4]\times \Gamma)\rightarrow L^2([0,T+\epsilon_1]\times \Gamma)}
\label{inequality-d1}
\eea
for some constant $C>0$.

\medskip
Step 3. We now estimate the two terms on the left hand side of the inequality (\ref{inequality-d1}).
First, by (\ref{condition-ad2}) and Lemma \ref{g-norm}, we have
\be \label{111}
|\langle \xi_1, \nu(x_1)\rangle \cdot |\bigl\langle g, g \bigr\rangle_{L^2(I(t_1,\epsilon_1/2)\times V(x_1))}| \thickapprox 1.
\ee
We next estimate $\langle \tilde{g}, g \rangle_{L^2(I(t_1,\epsilon_1/2)\times V(x_1))}$.
In the coordinate $x= F(y)$, by Lemma \ref{lem-g-ghat1},
we have
\be \label{inequality-d2}
\langle \hat{\tilde{g}}, \hat{g} \rangle_{L^2(I(t_1,\epsilon_1/2)\times U(x_1))} = \langle \hat{\tilde{g}}_*, \hat{g}_* \rangle_{L^2(I(t_1,\epsilon_1/2)\times U(x_1))} + O(\f{1}{\sqrt{\lambda}})
\ee
where $\hat{g}_*$ and $\hat{\tilde{g}}_*$ are the auxiliary beams associated with the beams $\hat{g}$ and $\hat{\tilde{g}}$ respectively. We write $\hat{g}_*= a(t_1) e^{i \lambda \hat{\tau}_*}$ and $\hat{\tilde{g}}_*=  \tilde{a}(t_1) e^{i \lambda \hat{\tilde{\tau}}_*}$ with
\beas
\hat{\tau}_* &=& \langle (-1, \f{\p F}{\p y}(0)^{\dag}\xi_1), (t-t_1, y)\rangle + (t-t_1, y)\hat{M}(t_1)(t-t_1, y)^{\dag}; \\
\hat{\tilde{\tau}}_* &=& \langle (-1, \f{\p F}{\p y}(\delta y)^{\dag} \tilde{\xi}_1), (t-\tilde{t}_1, y-\delta y)\rangle + (t-\tilde{t}_1, y-\delta y)\hat{\tilde{M}}(t_1)(t-\tilde{t}_1, y- \delta y)^{\dag}.
\eeas
By Lemma \ref{inequality2c} in Appendix \ref{appendix-beam-interaction}, we have
\be \label{inequality-d3}
|\langle \hat{\tilde{g}}_*, \hat{g}_* \rangle_{L^2(I(t_1,\epsilon_1/2)\times U(x_1))}| \ \lesssim e^{-c_0 \lambda |\delta z|}
\ee
where $c_0$ is a positive constant depending only on $\|c\|_{C^3}+ \|\tilde{c}\|_{C^3}$
and $|\delta z|= |t_1-\tilde{t}_1|^2 + |\delta y|^2+ |\f{\p F}{\p y}(\delta y)^{\dag} \tilde{\xi}_1- \f{\p F}{\p y}(0)^{\dag}\xi_1|^2$.
It follows from (\ref{inequality-d2}) and (\ref{inequality-d3}) that
\be \label{112}
|\langle \tilde{g}, g \rangle_{L^2(I(t_1,\epsilon_1/2)\times V(x_1))}| \ \lesssim e^{-c_0 \lambda |\delta z|}+ O(\f{1}{\sqrt{\lambda}}).
\ee

\medskip

Step 4. Combining (\ref{inequality-d1}), (\ref{111}) and (\ref{112}), we see that
$$
  e^{-c_0 \lambda |\delta z|} \gtrsim C_1- C_2 \|\Lambda_{\tilde{c}}- \Lambda_{c}\|_{H_0^1([0, 3\epsilon_1/4]\times \Gamma)\rightarrow L^2([0,T+\epsilon_1]\times \Gamma)} - C_3 \f{1}{\sqrt{\lambda}}
$$
for some positive constants $C_1$, $C_2$ and $C_3$ which are independent of $(x_0, \xi_0)$. By letting $\lambda \to \infty$, we conclude
that $\delta z= 0$ if
$$
\|\Lambda_{\tilde{c}}- \Lambda_{c}\|_{H_0^1([0, 3\epsilon_1/4]\times \Gamma)\rightarrow L^2([0,T+\epsilon_1]\times \Gamma)}< \f{C_1}{C_2}.
$$
Set $\delta=\f{C_1}{C_2}$.
From $\delta z= 0$ it follows that $t_1= \tilde{t}_1$, $\delta y=0$, and $\f{\p F}{\p y}(0)^{\dag} \tilde{\xi}_1- \f{\p F}{\p y}(0)^{\dag}\xi_1=0$. It remains to show that
$\tilde{\xi}_1=\xi_1$. Indeed, $\f{\p F}{\p y}(0)^{\dag} \tilde{\xi}_1- \f{\p F}{\p y}(0)^{\dag}\xi_1=0$ implies that the tangential component of $\tilde{\xi}_1$ and $\xi_1$ are equal. Besides, $|\xi_1| =|\tilde{\xi}_1|$. These together with (\ref{condition-ad2}) yield that $\tilde{\xi}_1=\xi_1$. This completes the proof of the theorem.

\section{Geodesic X-ray transform with caustics}

\subsection{Local properties of the normal operator $\mathfrak{N}$} \label{subsec-operator-germ}

In this subsection, we present some results about the local properties of the normal operator $\mathfrak{N}$ (see (\ref{normal-operator})).

From now on, we fix $x_*\in \Omega$. We first decompose $\mathfrak{N}$ locally
into two parts based on the separation of singularities of its Schwartz kernel.  Note that the map
$\phi(x_*, \cdot): \R^d \to \R^d$ is a diffeomorphism in a neighborhood of the origin.
In fact, we can check that
$\f{\p \phi(x_*, \cdot)}{\p \xi}(0)= c(x_*)\cdot Id$.
Similar to the proof of existence of uniformly normal neighborhood in Riemannian manifold \cite{Lee},
we can find $\epsilon_2>0$ and a neighborhood of $x_*$, say $\tilde{U}(x_*)\subset \R^d$,
such that
\be
\phi(x, \cdot)|_{|\xi|< 2 \epsilon_2} \,\, \mbox{is a diffeomorphism for any} \,\,x\in \tilde{U}(x_*).
\ee

Let $\chi_{\epsilon_2} \in C_0^{\infty}(\R)$ be such that $\chi(t)=1$ for $|t|< \epsilon_2$ and
$\chi(t)=0$ for $|t|>2\epsilon_2$. We then define
\bea
\mathfrak{N_1}f(x)&=& \int_{T^*_{x}\Omega}W(x, \xi)f(\phi(x, \xi))\chi_{\epsilon_2}(|\xi|)\,d\sigma_x(\xi),\\
\mathfrak{N_2}f(x)&=& \int_{T^*_{x}\Omega}W(x, \xi)f(\phi(x, \xi))(1-\chi_{\epsilon_2}(|\xi|))\,d\sigma_x(\xi).
\eea
Note that for any $f$ supported in $\Omega$, $f(\phi(x, \xi))=0$ for all $|\xi| > T$. Thus we have
$$
\mathfrak{N_2}f(x)= \int_{\xi\in T^*_{x}\Omega, \,\, \epsilon_2 < |\xi|< T} W(x, \xi)f(\phi(x, \xi))(1-\chi_{\epsilon_2}(|\xi|))\,d\sigma_x(\xi).
$$

It is clear that $\mathfrak{N} f= \mathfrak{N_1} f+\mathfrak{N_2} f$.
This gives the promised decomposition of $\mathfrak{N}$. We next study $\mathfrak{N_1}$ and $\mathfrak{N_2}$
separately.

\begin{lem} \label{prop-old1}
$\mathfrak{N_1}$ is an elliptic $\Psi$DO of order $-1$ from $C_0^{\infty}(\tilde{U}(x_*), \R^{2d})$ to $\mathcal{D}'(\tilde{U}(x_*), \R^{2d})$
 with principle symbol
\be
\sigma_p(\mathfrak{N}_1)(x ,\xi)=
2\pi \cdot \int_{S^*_x\Omega} \delta(\langle \xi, \theta \rangle) \Phi^\dag(x,\theta)\cdot \Phi(x,\theta)\,d\sigma_x(\theta). \nonumber
\ee
\end{lem}
{\bf{Proof}}. See \cite{SU08AJM} or \cite{SUAPDE}.

\bigskip

We now proceed to study the operator $\mathfrak{N_2}$ whose property is determined by the Lagrangian map
$\phi(x_*, \cdot)$.
We shall study the operator germ $\mathfrak{N}_{2, \xi_*}$ for each $\xi \in T_{x_*}^*\R^d$.
We first consider the case when
$\xi_*$ is not a caustic vector, i.e. $\xi_*$ is a regular vector.

\begin{lem} \label{lem-x-ray5}
Let $\xi_* \in S^*_{x_*}\R^d$ be a regular vector, then there exists a neighborhood $U(x_*)$ of $x_*$ and a neighborhood
$B(x_*, \xi_*)$ of $(x_*, \xi_*)$ such that for any
$\chi\in C^{\infty}_0(B(x_*, \xi_*))$ the following operator
\be
\mathfrak{N}_{2, \xi_*} f(x)= \int_{T^*_{x}\Omega}W(x, \xi)f(\phi(x, \xi))(1- \chi_{\epsilon_2}(|\xi|))\cdot \chi(x, \xi)\,d\sigma_x(\xi) \nonumber
\ee
is a smoothing operator from $\mathcal{E}'(\Omega, \R^{2d})$ into $C^{\infty}(\overline{U(x_*)}, \R^{2d})$.
\end{lem}

{\bf{Proof}}. Since $\xi_* \in S^*_{x_*}\R^d$ is regular,
there exist a neighborhood $V(x_*)$ of $x_*$ in $\R^d$ and a neighborhood $B(x_*, \xi_*)$ of $(x_*, \xi_*)$ in $\R^{2d}$ of the form
$B(x_*, \xi_*)=V(x_*) \times B_0(\xi_*)$ for some open set $B_0(\xi_*)$ in $\R^d$
such that the map
$\phi(x, \cdot)$ is a diffeomorphism between $B_0(\xi_*)$ and its image for all $x\in V(x_*)$. We denote the inverse of the map
$\phi(x, \cdot)$ by $\phi^{-1}(x,\cdot)$. By a change of coordinate $\xi= \phi^{-1}(x, y)$ and use some cut-off function,
we can write $\mathfrak{N}_{2, \xi_*}$ in the following form
$$
\mathfrak{N}_{2, \xi_*} f(x)= \int_{\Omega}K(x, y) f(y)\,dy, f\in \mathcal{E}'(\Omega, \R^{2d})
$$
for some smooth function $K$ in $\Omega \times \Omega$. The Lemma follows immediately.

\medskip

We next consider the case when $\xi_*$ is a fold vector.
We have the following slightly modified result from \cite{SUAPDE}.

\begin{lem} \label{lem-modified-n2}
Let $\xi_{*}$ be a fold vector of the map $\phi(x_*,\cdot)$.
Then there exists a small neighborhood $U(x_*)$ of $x_*$
and a small neighborhood $B(x_*, \xi_*)$ of $(x_*, \xi_*)$ in $\R^{2d}$
such that for any $\chi \in  C^{\infty}_{0}(B(x_*, \xi_*))$, the operator
$\mathfrak{N}_{2, \xi_*}:  \mathcal{E}'(\Omega, \R^{2d}) \to  \mathcal{D}'(U(x_*), \R^{2d})$ defined by
\be
\mathfrak{N}_{2, \xi_*} f(x)= \int_{T^*_{x}\Omega}W(x, \xi)f(\phi(x, \xi))(1- \chi_{\epsilon_2}(|\xi|))\cdot \chi(x, \xi) \,d\sigma_x(\xi),
\quad f\in \mathcal{E}'(\Omega, \R^{2d})
\ee
is an FIO of order $-\f{d}{2}$ whose associated canonical relation is compactly supported in the following set
\bea \label{canonical-relation}
\Big\{(x, \xi, y, \eta);  && x\in U(x_*), y=\phi(x, \omega), (x, \omega) \in B(x_*, \xi_*), \, \, \det d_{\omega}\phi(x, \omega)=0,  \nonumber\\
&& \xi= -\eta_i\f{\p \phi^i(x, \omega)}{\p x}, \, \eta \in \mbox{Coker} \,(d_{\omega}\phi(x, \omega)).\Big\}
\eea
\end{lem}

\subsection{Singularities of the map $\phi(x, \cdot)$}
In this subsection, we present some properties about the map $\phi(x, \cdot)$ which is equivalent to the exponential map in Riemannian manifold.

By the classification result for Lagrangian maps
(see \cite{A} and \cite{AZV} for detail),
there are only a finite number of
stable and simple singular Lagrangian map germs in dimensions between three and five and they are generic.
In three dimensions, there are four types: fold, cusp, swallow-tail and D4. The others are unstable and can be removed by using arbitrarily small perturbations. We define
 \beas
 \mathcal{K}(x)&=& \{\xi\in T^*_{x} \R^d : \mbox{the map germ $\phi(x, \cdot)$ at $\xi$ is singular}\}; \\
 \mathcal{K}_1(x)&=& \{\xi\in T^*_{x}\R^d : \mbox{the map germ $\phi(x, \cdot)$ at $\xi$ has singularity of fold type}\}; \\
 \mathcal{K}_2(x)&=& \{\xi\in T^*_{x}\R^d : \mbox{the map germ $\phi(x, \cdot)$ at $\xi$ has singularity of cusp type}\}; \\
 \mathcal{K}_3(x)&=& \{\xi\in T^*_{x}\R^d : \mbox{the map germ $\phi(x, \cdot)$ at $\xi$ has simple and stable singularities of types} \\
                 &&  \mbox{\quad other than fold and cusp}\}.
 \eeas


We say that the map $\phi(x, \cdot)$ is in a general position (or generic) if the map germ $\phi(x, \cdot)$ is simple and stable at all caustic vectors in $\mathcal{K}(x)$.
By the classification result of Lagrangian maps, see for instance \cite{A},
the following result holds for the set $\mathcal{K}(x)$.

\begin{prop}
Assume that the map $\phi(x, \cdot)$ is in a general position, then the sets $\mathcal{K}_1(x)$ and $\mathcal{K}_2(x)$ are smooth manifolds of dimensions $d-1$ and $d-2$, respectively. The set $\mathcal{K}_3(x)$ is a union of smooth manifolds of dimensions not greater than $d-3$. Especially, for $d=3$, the sets $\mathcal{K}_1(x)$, $\mathcal{K}_2(x)$ and $\mathcal{K}_3(x)$ consists of smooth surfaces, smooth curves and isolated points,  respectively.
\end{prop}

In the case when the map $\phi(x, \cdot)$ is not in a general position, it is known that
$\bigcup_{j=1}^3 \mathcal{K}_j(x)$ is open and dense in $\mathcal{K}(x)$.

Recall that the map $\beta$ is defined by (\ref{beta}). Denote $\mathcal{S}_j(x)=\beta (\mathcal{K}_j(x))$ for j=1,2,3.
We conclude that the following result holds.

\begin{lem} \label{caustic-class}
Assume that the map $\phi(x, \cdot)$ is in general position, then the sets $\mathcal{S}_1(x)$, $\mathcal{S}_2(x)$  and $\mathcal{S}_3(x)$ are of finite $d-1$, $d-2$ and $d-3$ dimensional Hausdorff measures, respectively.
Especially, for $d=3$, the set $\mathcal{S}_2(x)$ is a curve (not necessarily smooth) of finite length in $S^*_{x} \R^3$ and $\mathcal{S}_3(x)$ consists of a finite number of points.
\end{lem}

It is clear that $\mathcal{S}_1(x)$ represents the set of geodesics passing through $x$ with fold caustics (and possibly other types of caustics as well). Observe that for each fold vector $\xi\in T_x^*\R^d$, there exists a neighborhood of $\xi$ in $\mathcal{K}_{1}(x)$, say $V(\xi)$, such that $V(\xi)$ is a smooth hyperplane of dimension $d-1$ in $T_x^*\R^d$ and is transversal to the ray $\{t\xi: t>0\}$. From this,
we can conclude that the following result holds.
\begin{lem} \label{caustic-class1}
The set $\mathcal{S}_1(x)$ is open in $S^*_{x} \R^d$.
\end{lem}

\subsection{Discussions on the concept of Fold-regular}  \label{subsec-discussion}

In this subsection, we discuss the concept ``fold-regular'' defined in Subsection \ref{subsec-fold-regular}. We show that fold-regular point is the natural object to study for geodesic X-ray transforms with caustics. We also derive some natural conditions for a point to be fold-regular.
We remark that one can define ``fold-regular metric'' in a similar way as for ``fold-regular velocity field'' in Subsection \ref{subsec-thm3}. Then ``fold-regular metric'' generalizes ``regular metric'' introduced in \cite{SU08JGA}.


We take the case $d=3$ for example, similar discussions also hold for the case $d>3$. Note that the set of geodesics passing through $x$ can be parameterized by the sphere $S_x^*\R^3$ and those with fold caustics by $\mathcal{S}_1(x)$, which is a subset with the same dimension. In comparison, a complete set of geodesics passing through $x$ has dimension at least one. 
By Lemma \ref{caustic-class1}, we see that in the generic case for a non-simple velocity field, the set
$\mathcal{S}_1(x)$ is open in the sphere $S^*_x\R^3$. Thus, in order to be able to select a complete set from $S^*_x\R^3$, it is necessary to use geodesics with fold caustics in the case when there is no complete subset in the set of geodesics without caustics.  However, we should also point out that there is not much need to use geodesics with other types of caustics since they correspond to a one dimensional subset on $S^*_x\R^3$. In the extreme case, when the set $\mathcal{S}_1(x)$ contains no complete subset, then there exists $\theta \in S^*_{x}\R^3$ such that all the geodesics represented by the set $\theta^{\perp} = \{\xi \in S^*_{x}\R^d, \xi \perp \theta\}$ have caustics other than the fold type. It is only in such extreme case that we need to use information from geodesics with cusp caustics and hence an investigation of properties of the operator germ $\mathfrak{N}_{\xi}$ for cusp vectors becomes necessary.


We now consider the problem when a fold vector is fold-regular.
In dimension $d \geq 3$, a sufficient condition for a fold vector $\xi_*$ to be fold-regular is the following
\be \label{graph-condition}
  d^2_{\xi} \phi(x_*, \xi_*)(N_{x_*}(\xi_*)\setminus 0 \times \cdot )|_{T_{\xi_*}S(x_*)} \quad \mbox{is of full rank},
\ee
where $N_{x_*}(\xi_*)$ denotes the kernel of $d_{\xi} \phi(x_*, \xi_* )$  and
$S(x_*)$ the set of all vectors $ \xi \in T^*_{x_*} \R^d$ such that $\det d_{\xi} \phi(x_*, \xi)=0$.
Indeed, in that case, it is shown in \cite{SUAPDE} that the canonical relation associated with the operator germ $\mathfrak{N}_{2, \xi_*}$ is locally a canonical graph and hence $\mathfrak{N}_{2, \xi_*}$ is bounded from $L^2(\Omega_{\epsilon_0}, \R^{2d})$ to
$H^{\f{d}{2}}(U( x_*), \R^{2d})$ for some neighborhood $U(x_*)$ of $x_*$. Note that for $d\geq 3$,
$H^{\f{d}{2}}(U(x_*), \R^{2d})$ is compactly embedded in $H^{1}(U(x_*), \R^{2d})$, so $\mathfrak{N}_{2, \xi_*}$ is compact from
$L^2(\Omega_{\epsilon_0}, \R^{2d})$ to $H^1(U( x_*), \R^{2d})$ and we can conclude that $\xi_*$ is fold-regular.

The set of fold-regular vectors contains more elements than those which satisfy the graph condition (\ref{graph-condition}). In fact, let $\mathcal{C}\subset T^*\Omega \times  T^*\Omega$ be the canonical relation associated
with the operator germ $\mathfrak{N}_{2, \xi_*}$ defined in Lemma \ref{lem-modified-n2}.
We can show that $\mathcal{C}$ is homogeneous and
$\mathcal{C} \subset (T^*\Omega \backslash 0) \times  (T^*\Omega \backslash 0)$. By the main result in \cite{GS-AJM-98},
$\mathfrak{N}_{2, \xi_*}$ is bounded from $L^2(\Omega_{\epsilon_0}, \R^{2d})$ to $H^{\f{d}{2}-\f{1}{3}}(U(x_*), \R^{2d})$
for some neighborhood $U(x_*)$ of $x_*$,
if the only singularity of
the projection of  $\mathcal{C}$
to its first or second component at the point associated with $(x_*, \xi_*)$ is fold or cusp. Since $H^{\f{d}{2}-\f{1}{3}}(U(x_*), \R^{2d})$ is compactly embedded in
$H^{1}(U(x_*), \R^{2d})$, we see that $\mathfrak{N}_{2, \xi_*}$ is compact from
$L^2(\Omega_{\epsilon_0}, \R^{2d})$ to $H^1(U(x_*), \R^{2d})$ and hence $\xi_*$ is fold-regular.

\begin{rmk}
In dimension $d=2$, the set of fold-regular vectors is generally empty. Indeed, for a fold vector $\xi_*$,
the operator germ $\mathfrak{N}_{2, \xi_*}$ is a FIO of order $-1$, and hence the best estimate is that it is bounded from
 $L^2(\Omega_{\epsilon_0}, \R^{2d})$ to $H^1(U(x_*), \R^{2d})$ for some neighborhood $U(x_*)$ of $x_*$.
\end{rmk}




We conclude this subsection with the following criterion for a fold-regular point.
\begin{lem}
If for any $\theta \in S_x^*\R^d$, there exists $\xi \in S_x^*\R^d$ such that $\xi \perp \theta$ and that either no caustics or   only fold caustics with condition (\ref{graph-condition}) satisfied exist along the ray $\{ \phi(x, t\xi): -T \leq t \leq T\}$, then the point $x$ is fold-regular.
\end{lem}

\subsection{Proof of Theorem \ref{thm2}}

We prove Theorem \ref{thm2} in this subsection. The proof can be divided into two major stages: in the first stage, we present some preliminaries and construct a cut-off function $\alpha \in C_0^{\infty}(S^*_{-}\Gamma)$ which selects a complete set of geodesics with only fold-regular caustics, see Lemma \ref{lem-cutoff}; in the second stage, we study the normal operator $\mathfrak{N}_{\alpha}=\mathfrak{I}_{\alpha}^{\dag}\mathfrak{I}_{\alpha}$, see Lemma \ref{lem-n1} and \ref{lem-n2}. Theorem \ref{thm2} is then a direct consequence of Lemma \ref{lem-n1} and  \ref{lem-n2}.

We now present some preliminaries that are necessary for the construction of $\alpha$.
Let $x_*$ be a fold-regular point with the compact subset $\mathcal{Z}(x_*)\subset S^*_{x_*}\R^d$ in Definition \ref{def-fold-regular-point}.
Denote
$\mathcal{C}_{\epsilon_2, T}\mathcal{Z}(x_*)= \{r\xi;  \xi\in \mathcal{Z}(x_*), r\in \R \quad \mbox{and} \,\,\, \epsilon_2\leq |r|\leq T\}$.
For each $\xi_* \in \mathcal{C}_{\epsilon_2, T}\mathcal{Z}(x_*)$, by Lemma \ref{lem-x-ray5} and Lemma \ref{lem-modified-n2}, there exist a neighborhood $U(x_*, \xi_*)$ of $x_*$ and a neighborhood $B(x_*, \xi_*)$ of $(x_*, \xi_*)$  such that for any
$\chi \in C^{\infty}_0(B(x_*, \xi_*))$ the following operator
\[
\mathfrak{N}_{2, \xi_*} f(x)= \int_{T^*_{x}\Omega}W(x, \xi)f(\phi(x, \xi))(1- \chi_{\epsilon_2}(|\xi|))\cdot \chi(x,\xi)\,d\sigma_x(\xi)
\]
is compact from $H^s(\Omega_{\epsilon_0}, \R^{2d})$ to $H^{s+1}(U(x_*, \xi_*), \R^{2d})$. Let $B_0(x_*, \xi_*)$ be another neighborhood of $(x_*, \xi_*)$ in $\R^{2d}$
such that $B_0(x_*, \xi_*) \Subset B(x_*, \xi_*)$.
Since $\mathcal{C}_{\epsilon_2, T}\mathcal{Z}(x_*)$ is compact, there exists a finite number of $\xi_*$'s in $\mathcal{C}_{\epsilon_2, T}\mathcal{Z}(x_*)$, say $\xi_1$,  $\xi_2$,...,  $\xi_M$, such that
\[
\mathcal{C}_{\epsilon_2, T} \mathcal{Z}(x_*) \subset \bigcup_{j=1}^{M} B_0(x_*, \xi_j).
\]
We can then find smooth functions
$\chi_1$, $\chi_2$,...,$\chi_M$ with $\mbox{supp}\chi_j \subset B(x_*, \xi_j)$ for each $j$ such that
$$
\sum_{j=1}^{M}\chi_j(x, \xi)= 1 \quad \mbox{for all } \,\, (x, \xi) \in \bigcup_{j=1}^{M} B_0(x_*, \xi_j).
$$

Denote by $\mathcal{A}_0$ be the greatest connected open symmetric subset in $\bigcup_{j=1}^{M} B_0(x_*, \xi_j)$
which contains $\mathcal{C}_{\epsilon_2, T} \mathcal{Z}(x_*) $. Here and after, we say that a set $\mathcal{B}$
in $\R^{2d}$ is symmetric if $(x, \xi)\in \mathcal{B}$ implies that $(x, -\xi)\in \mathcal{B}$.
Define
$$
\mathcal{A}_{\epsilon} =\{(x, \xi) \in \R^{2d}: |x-x_*|\leq \epsilon, \,\, \epsilon_2 \leq |\xi| \leq T\}
$$
for each $\epsilon>0$. It is clear that $\mathcal{A}_{\epsilon}$ is compact in $\R^{2d}$,
so is the set $\mathcal{A}_{\epsilon}\backslash \mathcal{A}_0$.

\begin{lem} \label{lem-cutoff}
There exist $\epsilon_3>0$ and $\alpha \in C_0^{\infty}(S^*_{-}\Gamma)$ such that the following two conditions are satisfied:
\bea
\alpha(x_0, \xi_0)&=& 1\quad \mbox{for all } \, (x_0, \xi_0)\in \tau\circ \beta(\mathcal{C}_{\epsilon_2, T}\mathcal{Z}(x_*)), \label{alpha1}\\
\alpha(x_0, \xi_0)&=& 0 \quad \mbox{for all } \, (x_0, \xi_0)\in \tau\circ\beta(\mathcal{A}_{\epsilon_3}\backslash \mathcal{A}_0)\label{alpha2}.
\eea
\end{lem}
{\bf{Proof}}.
Note that both $\beta$ and $\tau$ are continuous. Since $\mathcal{C}_{\epsilon_2, T}\mathcal{Z}(x_*)$
and $\mathcal{A}_{\epsilon}\backslash \mathcal{A}_0$ are compact, so are the sets
$\tau\circ\beta(\mathcal{C}_{\epsilon_2, T}\mathcal{Z}(x_*))$ and $\tau\circ\beta(\mathcal{A}_{\epsilon}\backslash \mathcal{A}_0)$. We claim that
there exists $\epsilon_3>0$ such that
$$
\tau\circ \beta(\mathcal{C}_{\epsilon_2, T}\mathcal{Z}(x_*)) \bigcap \tau\circ\beta(\mathcal{A}_{\epsilon}\backslash \mathcal{A}_0) = \emptyset
$$
for all $\epsilon\leq \epsilon_3$. Indeed, assume the contrary, then
$$
\tau\circ \beta(\mathcal{C}_{\epsilon_2, T}\mathcal{Z}(x_*)) \bigcap \tau\circ\beta(\mathcal{A}_{\epsilon}\backslash \mathcal{A}_0) \neq \emptyset
$$
for all $\epsilon>0$. Note that
the collection of compact subsets
$\tau\circ \beta(\mathcal{C}_{\epsilon_2, T}\mathcal{Z}(x_*)) \bigcap \tau\circ\beta(\mathcal{A}_{\epsilon}\backslash \mathcal{A}_0)$ is decreasing with respect to $\epsilon$, so it satisfies the finite intersection property and we can thus conclude that
$$
\tau\circ \beta(\mathcal{C}_{\epsilon_2, T}\mathcal{Z}(x_*)) \bigcap_{\epsilon>0}
\tau\circ\beta(\mathcal{A}_{\epsilon}\backslash \mathcal{A}_0) \neq \emptyset.
$$
But on the other hand, we can check that
\beas
\bigcap_{\epsilon>0} \tau\circ\beta(\mathcal{A}_{\epsilon}\backslash \mathcal{A}_0)
&=& \tau( (\mathcal{A}_{\epsilon}\backslash \mathcal{A}_0) \bigcap  S^*_{x_*}\R^d)  \\
\tau\circ \beta(\mathcal{C}_{\epsilon_2, T}\mathcal{Z}(x_*))
&=& \tau(\mathcal{C}_{\epsilon_2, T}\mathcal{Z}(x_*) \bigcap S^*_{x_*}\R^d).
\eeas
Using the fact that $\tau$ is injective on $S^*_{x_*}\R^d$ and $\mathcal{C}_{\epsilon_2, T} \mathcal{Z}(x_*) \subset \mathcal{A}_0$, we obtain
$$
\tau( (\mathcal{A}_{\epsilon}\backslash \mathcal{A}_0) \bigcap  S^*_{x_*}\R^d) \bigcap
\tau(\mathcal{C}_{\epsilon_2, T}\mathcal{Z}(x_*) \bigcap S^*_{x_*}\R^d) = \emptyset.
$$
Thus,
$$
\tau\circ \beta(\mathcal{C}_{\epsilon_2, T}\mathcal{Z}(x_*)) \bigcap_{\epsilon>0}
\tau\circ\beta(\mathcal{A}_{\epsilon}\backslash \mathcal{A}_0)  =\emptyset.
$$
This contradiction completes the proof of our claim.

Now, we have
$$
\tau\circ \beta(\mathcal{C}_{\epsilon_2, T}\mathcal{Z}(x_*)) \bigcap
\tau\circ\beta(\mathcal{A}_{\epsilon_3}\backslash \mathcal{A}_0)  = \emptyset.
$$
By decreasing $\epsilon_3$ if necessary, we may assume that
$$
\{x: |x-x_*| \leq \epsilon_3\} \subset  \pi(\mathcal{A}_0).
$$
Since both the sets $\tau\circ \beta(\mathcal{C}_{\epsilon_2, T}\mathcal{Z}(x_*))$ and
$\tau\circ\beta(\mathcal{A}_{\epsilon_3}\backslash \mathcal{A}_0)$
are compact in $S^*_{-}\Gamma$, we can find
$\alpha\in C_0^{\infty}(S^*_{-}\Gamma)$ as desired. This concludes the proof of the lemma.

\medskip

The construction of $\alpha$ above completes the first stage of the proof of Theorem \ref{thm2}, we are now at the second stage.
We define the truncated geodesic X-ray transform $\mathfrak{I}_{\alpha}f$ as in (\ref{formula-I-alpha}) or (\ref{formula-I-alpha-1}).
By replacing the weight $\Phi$ with the new one $\alpha^{\sharp}\cdot \Phi$, we obtain
$\mathfrak{N}_{\alpha}$, $\mathfrak{N}_{1, \alpha}$ and $\mathfrak{N}_{2, \alpha}$ from the corresponding formulas of
$\mathfrak{N}$, $\mathfrak{N}_{1}$ and $\mathfrak{N}_{2}$.
It is clear that Lemma \ref{lem-x-ray5}, \ref{lem-modified-n2} still hold with the new weight.

\begin{lem}\label{lem-n1}
There exist a neighborhood $U(x_*)$ of $x_*$ and
a smoothing operator $\mathfrak{R}$ from $\mathcal{E}'(\Omega, \R^{2d})$ into $C^{\infty}(\overline{U(x_*)}, \R^{2d})$,
such that for for any $s\in \R$ and any neighborhood $U_0(x_*)$ of $x_*$ with $U_0(x_*) \Subset U(x_*)$, the following estimate holds
\be
\|f\|_{H^s(U_0(x_*), \R^{2d})} \lesssim \|\mathfrak{N}_{1, \alpha}f \|_{H^{s+1}(U(x_*), \R^{2d})} + \|\mathfrak{R}f\|_{H^{s}(\Omega, \R^{2d})}.
\ee
\end{lem}

{\bf{Proof}}. We first show that $\mathfrak{N}_{1, \alpha}$ is an elliptic $\Psi$DO .
Indeed, as in Lemma \ref{prop-old1}, $\mathfrak{N}_{1, \alpha}$ is a
$\Psi$DO of order $-1$ from $C_0^{\infty}(\tilde{U}(x_*), \R^{2d})$ to $\mathcal{D}'(\tilde{U}(x_*), \R^{2d})$
with principle symbol
\be
\sigma_p(\mathfrak{N}_1)(x ,\xi)=
2\pi \cdot \int_{S^*_x\Omega} \delta(\langle \xi, \theta \rangle) |\alpha^{\sharp}(x_*,\theta)|^2\Phi^\dag(x,\theta)\cdot \Phi(x,\theta)\,d\sigma_x(\theta). \nonumber
\ee

By the construction of $\alpha$, for any $\xi\in S^*_{x_{*}}\R^d$,  we have $\alpha^{\sharp}(x_*,\theta)=1$ for some
$\theta \in S^*_{x_*}\R^d$ with $\theta \perp \xi$. Thus
$$
\sigma_p(\mathfrak{N_{1, \alpha}})(x_* ,\xi)=
2\pi \cdot \int_{\theta \in S^*_{x_*}\R^d, \, \theta \perp \xi} |\alpha^{\sharp}(x_*,\theta)|^2\Phi^\dag(x_*,\theta)\cdot \Phi(x_*,\theta)\,d\sigma_{x_*}(\theta) >0
$$
in the sense of symmetric positive definite matrix. By continuity, we can find a neighborhood $U(x_*)\subset \tilde{U}(x_*)$ of $x_*$ such that
$\sigma_p(\mathfrak{N}_{1, \alpha})(x , \xi)> 0$ for all $x\in U(x_*)$ and $\xi \in S^*_{x}\R^d$. Thus we can conclude that
$\mathfrak{N}_{1, \alpha}$ is an
elliptic $\Psi$DO of order $-1$ from $C_0^{\infty}(U(x_*), \R^{2d})$ to $\mathcal{D}'(U(x_*), \R^{2d})$. By standard argument, we can conclude that
Lemma \ref{lem-n1} holds.

\bigskip

We now study the operator $\mathfrak{N}_{2, \alpha}$.

\begin{lem} \label{lem-n2}
There exists a small neighborhood of $x_*$, say $U(x_*)$, such that the operator
$\mathfrak{N}_{2, \alpha}: \mathcal{E}'(\Omega, \R^{2d}) \to \mathcal{D}'(U(x_*), \R^{2d})$
is compact from $H^s(\Omega_{\epsilon_0}, \R^{2d})$ to $H^{s+1}(U(x_*), \R^{2d})$.
\end{lem}

{\bf{Proof}}. Recall that $\mathfrak{N}_{2, \alpha}$ has the following representation
\[
\mathfrak{N}_{2, \alpha} f(x)= \int_{T^*_{x}\Omega}W_{\alpha}(x, \xi)f(\phi(x, \xi))(1- \chi_{\epsilon_2}(|\xi|))\cdot \,d\sigma_x(\xi),
\]
where
\beas
W_{\alpha}(x, \xi) &=&  \f{1}{|\xi|^{d-1}}
  |\alpha\circ \tau\circ \beta(x, \xi)|^2\Phi^\dag \circ \beta(x, \xi) \cdot \Phi\circ \beta\circ \mathcal{H}(x, \xi)\\
 && + \f{1}{|\xi|^{d-1}}|\alpha\circ \tau\circ \beta(x, -\xi)|^2\Phi^\dag \circ \beta(x, -\xi) \cdot \Phi\circ \beta\circ \mathcal{H}^{-1}(x, -\xi).
\eeas
By (\ref{alpha2}) and the fact that $\mathcal{A}_0$ is symmetric,  we see that $\mbox{supp }W_{\alpha} \subset \mathcal{A}_0$ for all $x$ with
$|x-x_*| \leq \epsilon_3$.

Now, let $\chi_j$'s be as in the first stage. Define $\mathfrak{N}_{2, j}: \mathcal{E}'(\Omega, \R^{2d}) \to\mathcal{D}'(U(x_*, \xi_j), \R^{2d})$ by
\[
\mathfrak{N}_{2, j} f(x)= \int_{T^*_{x}\Omega}W_{\alpha}(x, \xi)f(\phi(x, \xi))(1- \chi_{\epsilon_2}(|\xi|))\cdot \chi_j(x,\xi)\,d\sigma_x(\xi).
\]

Let $U(x_*)= \bigcap_{j=1}^{M} (U(x_*, \xi_j)) \bigcap \{x: |x-x_*| < \epsilon_3\}$.
Then $U(x_*)$ is a neighborhood of $x_*$ and each $\mathfrak{N}_{2, j}$ is compact from $H^s(\Omega_{\epsilon_0}, \R^{2d})$ into
$H^{s+1}(U(x_*), \R^{2d})$.

We claim that $\mathfrak{N}_{2, \alpha}= \sum_{j=1}^M\mathfrak{N}_{2, j}$ when both sides are viewed as operators
from $\mathcal{E}'(\Omega, \R^{2d})$ to $\mathcal{D}'(U(x_*), \R^{2d})$. Indeed, for any $f\in C_0^{\infty}(\Omega, \R^{2d})$,
since $\sum_{j=1}^M\chi_j=1$ on $\mathcal{A}_0$ and $\mbox{supp }W_{\alpha} \subset \mathcal{A}_0$, we have
$$
W_{\alpha}(x, \xi)f(\phi(x, \xi))(1- \chi_{\epsilon_2}(|\xi|))
=W_{\alpha}(x, \xi)f(\phi(x, \xi))(1- \chi_{\epsilon_2}(|\xi|))\cdot \bigl(\sum_{j=1}^M\chi_j(x,\xi)\bigr)
$$
for all $x\in U(x_*)$. Thus $\mathfrak{N}_{2, \alpha}f= \sum_{j=1}^M\mathfrak{N}_{2, j}f$ and the claim follows.
This completes the proof of the lemma.

\medskip

Finally, note that $\mathfrak{N}_{\alpha}=\mathfrak{N}_{1, \alpha}+ \mathfrak{N}_{2, \alpha}$.
Theorem \ref{thm2} follows from Lemma \ref{lem-n1} and \ref{lem-n2}.

\section{Sensitivity analysis of recovering the velocity field from the DDtN map}
In this section, we prove Theorem \ref{thm3} on the
sensitivity of the inverse problem of recovering the velocity field from the DDtN map.
We first present a lemma which is a direct consequence of Theorem \ref{thm1} and Lemma \ref{lem-scattering-relation}.
\begin{lem} \label{lem-nonlinear}
Let $c$ and $\tilde{c}$ be two velocity field in
$\mathfrak{A}(\epsilon_0, \Omega, M_0, T)$, and let
$f$ be as in (\ref{f}).
Then there exists $\delta>0$ such that
if $\|\Lambda_{\tilde{c}}- \Lambda_{c}\|_{H_0^1[0, 3\epsilon_1/4]\times \Gamma \rightarrow L^2([0,T+\epsilon_1]\times \Gamma)} \leq \delta$, then
\be
\|\mathfrak{I}f\|_{L^{\infty}(S_{-}^*\Gamma, \R^{2d})} \leq C \|f\|_{C^1(\Omega, \R^{2d})}^2
\ee
for constant $C>0$ depending $M_0$.

\end{lem}

{\bf{Proof of Theorem \ref{thm3}}}. The proof is divided into the following six steps.

Step 1. Since $c$ is fold-regular, for each $x \in \overline{\Omega}_{\epsilon_0}$,
by Theorem \ref{thm2}, there exist a neighborhood $U(x)$ of $x$, a smooth cut-off function $\alpha \in C^{\infty}_0(S^*_{-}\Gamma)$
and a smoothing operator $\mathfrak{R}$ such that for any $U_0(x) \Subset U(x)$
the following estimate holds
\be
\|f\|_{L^{2}(U_0(x), \R^{2d})}
 \lesssim  \|\mathfrak{N}_{\alpha}f\|_{H^{1}(U(x), \R^{2d})} + \|\mathfrak{N}_{2, \alpha}f\|_{H^{1}(\Omega, \R^{2d})}+
 \|\mathfrak{R}f\|_{H^{1}(\Omega, \R^{2d})}
\ee
for all $f\in L^2(\Omega_{\epsilon_0}, \R^{2d})$.
Moreover, both $\mathfrak{N}_{2, \alpha}$ and $\mathfrak{R}$ are compact from $L^2(\Omega_{\epsilon_0}, \R^{2d})$ to $H^{1}(U(x), \R^{2d})$.

We now fix a neighborhood $U_0(x)\Subset U(x)$ of $x$ for each $x$. Note that $\overline{\Omega}_{\epsilon_0}$ is compact,
there exists a finite number of points, say $x_1$, $x_2$, ... $x_M$ such that
$\overline{\Omega}_{\epsilon_0} \subset \bigcup_{j=1}^{M}U_0(x_j)$.
Let $\mathfrak{N}_{\alpha_j}$ be the operator associated with each point $x_j$.

\medskip

Step 2. Denote by $H$ the Hilbert space $\prod_{j=1}^{M}H^{1}(U(x_j), \R^{2d})$. We consider the following three operators
\beas
T f &=& (\mathfrak{N}_{\alpha_1}f, \mathfrak{N}_{\alpha_2}f,..., \mathfrak{N}_{\alpha_M}f), \\
T_1f &=& (\mathfrak{N}_{2, \alpha_1}f, \mathfrak{N}_{2, \alpha_2}f,..., \mathfrak{N}_{2, \alpha_M}f),\\
T_2 f &=& (\mathfrak{R}_{\alpha_1}f, \mathfrak{R}_{\alpha_2}f,..., \mathfrak{R}_{\alpha_M}f).
\eeas
It is clear that all three operators are bounded from $L^2(\Omega_{\epsilon_0})$ to $H$.
Moreover, $T_1$ and $T_2$ are also compact and the following estimate holds
\be \label{inequality1}
\|f\|_{L^{2}(\Omega, \R^{2d})} \lesssim \|Tf\|_{H} + \|T_1f\|_{H} + \|T_2f\|_{H}.
\ee

Step 3. Let $\mathfrak{L}_0\subset L^2(\Omega_{\epsilon_0}, \R^{2d})$ be the kernel of $T$. We claim that
$\mathfrak{L}_0\subset L^2(\Omega_{\epsilon_0}, \R^{2d})$ is of finite dimension.
We prove by contradiction. Assume the contrary, then there exists an infinity number of orthogonal vectors in
$\mathfrak{L}_0\subset L^2(\Omega_{\epsilon_0}, \R^{2d})$,
say, $e_1$, $e_2$, ..., such that $\|e_j\|_{L^2(\Omega_{\epsilon_0}, \R^{2d})}=1$ and $Te_j=0$ for all $j\in \mathbb{N}$.
Since the sequence $\{e_j\}_{j=1}^{\infty}$ is bounded in $L^2(\Omega_{\epsilon_0}, \R^{2d})$ and the operators $T_1$ and $T_2$ are compact,
we can find a subsequence, still denoted by $\{e_j\}_{j=1}^{\infty}$,
such that both the sequences $\{T_1e_j\}_{j=1}^{\infty}$ and $\{T_2e_j\}_{j=1}^{\infty}$ are Cauchy in $H$.
By applying Inequality (\ref{inequality1}) to the vectors $e_i-e_j$ and recall that $T(e_i-e_j)=0$, we conclude that
the sequence $\{e_j\}_{j=1}^{\infty}$ is also Cauchy in $L^2(\Omega_{\epsilon_0}, \R^{2d})$. This contradicts to the fact that $\|e_i-e_j\|_{L^{2}(\Omega_{\epsilon_0}, \R^{2d})} >1$ for all $i\neq j$. This contradiction proves the claim.

\medskip

Step 4. Denote by $\mathfrak{L}_0^{\perp}$ the orthogonal space to $\mathfrak{L}_0$ in $L^2(\Omega_{\epsilon_0}, \R^{2d})$.
We claim that
\be \label{inequality2}
\|f\|_{L^{2}(\Omega_{\epsilon_0}, \R^{2d})} \lesssim \|Tf\|_{H} \quad \mbox{for all  }\,f\in \mathfrak{L}_0^{\perp}.
\ee
Indeed, assume the contrary, there exists a sequence $\{f_n\}_{n=1}^{\infty} \subset  \mathfrak{L}^{\perp}$ such that
$\|f_n\|_{L^2(\Omega_{\epsilon_0}, \R^{2d})} =1$ and  $\|Tf_n\|_{H}\leq \f{1}{n}$ for all $n$.
By the same argument as in Step 3, we can find
a subsequence, still denoted by $\{e_j\}_{j=1}^{\infty}$, such that both the sequences $\{T_1e_j\}_{j=1}^{\infty}$
and $\{T_2e_j\}_{j=1}^{\infty}$ are Cauchy in $H$.
By Inequality (\ref{inequality1}) and the fact that $\|Tf_n\|_{H}\leq \f{1}{n}$,
we can conclude that $\{f_n\}_{n=1}^{\infty}$ is also Cauchy in
$L^2(\Omega_{\epsilon_0}, \R^{2d})$. Let $f_0= \lim_{n\to \infty}f_n$,
then $\|Tf_0\|_{H} = \lim_{{n\to \infty}} \|Tf_n\|_{H}=0$. This implies that $f_0 \in \mathfrak{L}_0$.
However, note that $\mathfrak{L}_0^{\perp}$ is closed, as the limit of a sequence of functions in $\mathfrak{L}_0^{\perp}$,
$f_0$ must belong to
$\mathfrak{L}_0^{\perp}$. Therefore, we see that $f_0=0$.
But this contradicts to the fact that
$\|f_0\|_{L^2(\Omega_{\epsilon_0}, \R^{2d})}= \lim_{{n\to \infty}} \|f_n\|_{L^2(\Omega_{\epsilon_0}, \R^{2d})}=1$.
The claim is proved.

\medskip

Step 5. From now on, let $f$ be as in (\ref{f}).
We claim that
$$
\|Tf\|_{H} \lesssim  \|f\|_{H^{\f{15}{2}}(\Omega, \R^{2d})}^{\f{2}{5}} \cdot \|f\|_{L^{2}(\Omega, \R^{2d})}.
$$
Indeed, for each $\mathfrak{I}_{\alpha_j}$, by
Lemma \ref{lem-nonlinear}, we have
$\|\mathfrak{I}_{\alpha_j}f\|_{L^{\infty}(S_{-}^*\Gamma)} \leq C_1 \|f\|_{C^1}^2.$
Apply $\mathfrak{I}_{\alpha_j}^{\dag}$ to both sides and use the fact that $\mathfrak{I}_{\alpha_j}^{\dag}$
is bounded from $L^2$ to $L^2$(see \cite{S99}), we obtain
\be \label{nonlinear1}
\|\mathfrak{N}_{\alpha_j}f\|_{L^{2}(\Omega, \R^{2d})} \lesssim \|f\|_{C^1(\Omega, \R^{2d})}^2.
\ee
Then,
\beas
\|\mathfrak{N}_{\alpha_j}f\|_{H^{1}(U(x_j), \R^{2d})}
& \lesssim &  \|\mathfrak{N}_{\alpha_j}f\|_{H^{3}(U(x_j), \R^{2d})}^{\f{1}{3}} \cdot  \|\mathfrak{N}_{\alpha}f\|_{L^{2}(\Omega, \R^{2d})}^{\f{2}{3}}
\qquad (\mbox{by interpolation inequality})\\
& \lesssim &  \|\mathfrak{N}_{\alpha_j}f\|_{H^{3}(U(x_j), \R^{2d})}^{\f{1}{3}} \cdot  \|f\|_{C^1(\Omega, \R^{2d})}^{\f{4}{3}}
\qquad (\mbox{by (\ref{nonlinear1})}) \\
& \lesssim &  \|f\|_{H^{3}(\Omega, \R^{2d})}^{\f{1}{3}} \cdot  \|f\|_{C^1(\Omega, \R^{2d})}^{\f{4}{3}} \qquad (\mbox{by Lemma \ref{lem-n1}, \ref{lem-n2}})\\
& \lesssim &  \|f\|_{H^{3}(\Omega, \R^{2d})}^{\f{1}{3}} \cdot  \|f\|_{H^{3}(\Omega, \R^{2d})}^{\f{4}{3}} \qquad (\mbox{by interpolation inequality})\\
& = &  \|f\|_{H^{3}(\Omega, \R^{2d})}^{\f{5}{3}} \\
& \lesssim &   \|f\|_{H^{\f{15}{2}}(\Omega, \R^{2d})}^{\f{2}{5}} \cdot \|f\|_{L^{2}(\Omega, \R^{2d})}. \qquad (\mbox{by interpolation inequality})
\eeas

It follows that
\be  \label{inequality3}
\|Tf\|_{H} = \sum_{j=1}^M \|\mathfrak{N}_{\alpha_j}f\|_{H^{1}(U(x_j), \R^{2d})}
 \lesssim  \|f\|_{H^{\f{15}{2}}(\Omega, \R^{2d})}^{\f{2}{5}} \cdot \|f\|_{L^{2}(\Omega, \R^{2d})}.
\ee
This finishes the proof of our claim.

\medskip

Step 6. Denote by $\mathfrak{L}$
the projection of $\mathfrak{L}_0$ from $L^{2}(\Omega, \R^{2d})$ to the space $L^{2}(\Omega, \R^{d})$ by taking the
last three components. Note that the first three components of $f$ are zero, see (\ref{f}). Thus the condition
$\nabla (\ln c^2- \ln \tilde{c}^2) \perp \mathfrak{L}$ implies that $f \in \mathfrak{L}_0^{\perp}$. Consequently,
Inequality (\ref{inequality2}) holds. Combining this with (\ref{inequality3}), we see that
\be
\|f\|_{L^{2}(\Omega, \R^{2d})} \lesssim \|f\|_{H^{\f{15}{2}}(\Omega, \R^{2d})}^{\f{2}{5}} \cdot \|f\|_{L^{2}(\Omega, \R^{2d})}. \nonumber
\ee
Therefore, we must have $f=0$ for $\|f\|_{H^{\f{15}{2}}(\Omega, \R^{2d})}^{\f{2}{5}}$ sufficiently small.
Finally, note that $\|f\|_{H^{\f{15}{2}}(\Omega, \R^{2d})} \lesssim \|c - \tilde{c}\|_{H^{\f{17}{2}}(\Omega)}$ and that both
$c$ and $\tilde{c}$ vanishes near the boundary, we conclude that $f=0$ implies $c=\tilde{c}$.
This completes the proof of the theorem.

\section*{Acknowledgements} We thank the anonymous referees for many useful suggestions
which have significantly improved the presentation of the paper.

\section{Appendix}
\subsection{Linearization of ODE system} \label{appendix-ode}
Given the following ODE system:
$$
\dot{y}=f(y), \q y=y_o,
$$
where $f\in C^1(\R^d, \R^d)$. We consider the perturbed the system
$$
\dot{y_\epsilon}=f_\epsilon(y_\epsilon), \q y_\epsilon=y_o,
$$
where $f_{\epsilon}=f+ \epsilon g$ with $g\in C^1(\R^d)$.
We formally write $y_\epsilon(t)=y(t)+ r(t)= y(t)+ \epsilon \phi(t)+r_1(t)$, where
$\phi(0)=r(0)=r_1(0)=0$. By substituting $y_\epsilon(t)=y(t)+ \epsilon \phi(t)+r_1(t)$ into the perturbed system, we can derive that
$\phi$ satisfies the following equation:
\be \label{perturb1}
\dot{\phi}(t)= \f{\p f}{\p y}(y(t))\cdot \phi(t) + g(y(t)), \q \phi(0)=0.
\ee
By Grownwall's inequality, we can show that $|r(t)| \leq C \epsilon$ and $|r_1(t)|\leq C\epsilon^2$, where $C$ is a constant depending on
$\|f\|_{C^2} + \|g\|_{C^1}$.

We solve equation (\ref{perturb1}) as follows. Let $A(t)= \f{\p f}{\p y}(y(t))$. Let $\Phi(t)$ and $\Psi(t)$ be the solution to the following ODE system:
\begin{align*}
\dot{\Phi}(t)&= -\Phi(t)A(t), \q \Phi(0)= Id; \\
\dot{\Psi}(t)&= A(t)\Psi(t), \q \Psi(0)= Id.
\end{align*}

A straightforward calculation shows that $\Phi(t)\Psi(t)\equiv \Phi(0)\Psi(0)=Id$. Moreover,
$$
\phi(t)= \Phi(t)^{-1} \int_0^t \Phi(s)g(y(s))\,ds.
$$

\subsection{Proof of Lemma \ref{lem-g-ghat1} and Lemma \ref{lem-g-sum}} \label{appendix-proof}
\textbf{Proof of \ref{lem-g-ghat1}}.  We only show (\ref{g-approximation1-1}), since (\ref{g-approximation2-2}) follows in a similar way.
For simplicity, denote $D=(3\epsilon_1/4, t_1+ \epsilon_1/2)\times U(x_1)$. We first show that
\be \label{g-approximation1}
\hat{g}(t, y, \lambda) = \hat{g}_*(t, y,  \lambda) + O(\f{1}{\sqrt{\lambda}}) \quad \mbox{in }\,\, L^2(D).
\ee
Indeed, by direct calculation,
\be
\hat{g}(t, y,  \lambda)-\hat{g}_*(t, y,  \lambda)= (a(t)-a(t_1))e^{i\lambda \hat{\tau}_*}
+ a(t) (e^{i\lambda \hat{\tau}}- e^{i\lambda \hat{\tau}_*}).
\ee
It suffices to show that
\beas
R_1&:=& \|(a(t)-a(t_1))e^{i\lambda \hat{\tau}_*}\|_{L^2(D)} \lesssim \f{1}{\sqrt{\lambda}},\\
R_2&:=& \|a(t) (e^{i\lambda \hat{\tau}}- e^{i\lambda \hat{\tau}_*})\|_{L^2(D)} \lesssim \f{1}{\sqrt{\lambda}}.
\eeas
We first estimate $R_1$. By Lemma 3.1 in \cite{BQ1}, we have $|a(t)| \thickapprox \lambda^{\f{d}{4}}$. By Equation
(\ref{transport}), we further derive that $|\dot{a}(t)| \thickapprox \lambda^{\f{d}{4}} $, thus
$$
a(t)-a(t_1)= \int_0^1\dot{a}(t_1+ s(t-t_1))\,ds (t-t_1)= O(\lambda^{\f{d}{4}})|t-t_1|.
$$
Therefore,
\begin{align*}
\|(a(t)-a(t_1))e^{i\lambda \hat{\tau}_*}\|^2_{L^2(D)}
\lesssim \int_{D} \lambda^{\f{d}{2}}(t-t_1)^2 e^{-\lambda (t-t_1, y)\Im{\hat{M}}(t_1)(t-t_1, y)^{\dag}}\, dtdy
 \lesssim \f{1}{\lambda}.
\end{align*}
This proves $R_1 \lesssim \f{1}{\sqrt{\lambda}}$.

We next estimate $R_2$. Write $\hat{\tau}=\hat{\tau}_* + \delta \hat{\tau}$, then
$\delta \hat{\tau} = O(|(t-t_1, y)|^3)$ and hence
$|1- e^{i\lambda \delta \hat{\tau}}| \lesssim \lambda  \cdot O(|(t-t_1, y)|^3)$.
It follows that
$$
R_2   \leq \int_{D} |a(t)e^{i\lambda \hat{\tau}_*}|^2\cdot |1- e^{i \lambda \delta \hat{\tau}}| \,dtdy
\lesssim  \int_{D} \lambda^{\f{d}{2}}\cdot \lambda \cdot |(t-t_1, y)|^3
e^{-2 \lambda (t-t_1, y)\Im{\hat{M}}(t_1)(t-t_1, y)^{\dag}}\, dtdy
\lesssim  \f{1}{\lambda}.
$$
This completes the proof of (\ref{g-approximation1}).

We now proceed to show (\ref{g-approximation1-1}).
By direct calculation,
\beas
\f{\p \hat{g}}{\p y}-\f{\p \hat{g}_*}{\p y}
&=& i \lambda\f{\p \hat{\tau}}{\p y} \cdot \hat{g}- i \lambda\f{\p \hat{\tau}_*}{\p y} \cdot \hat{g}_* \\
&=& i \lambda(\f{\p \hat{\tau}}{\p y}- \f{\p \hat{\tau}_*}{\p y})\cdot \hat{g} + i \lambda\f{\p \hat{\tau}_*}{\p y} \cdot(\hat{g}-\hat{g}_*)
\eeas
One can check that $\f{\p \hat{\tau}}{\p y}- \f{\p \hat{\tau}_*}{\p y}=O|(t-t_1, y)|^2$, then a similar argument as used in the estimate of $R_1$ above shows that
$$
\|\lambda(\f{\p \hat{\tau}}{\p y}- \f{\p \hat{\tau}_*}{\p y})\cdot \hat{g}\|^2_{L^2(D)}
\lesssim  1.
$$
Besides, (\ref{g-approximation1}) implies that
\[
\|\lambda\f{\p \hat{\tau}_*}{\p y} \cdot(\hat{g}-\hat{g}_*)\|^2_{L^2(D)}
\lesssim \lambda.
\]
Combining these two estimates together, we conclude that
$$
\|\f{\p \hat{g}}{\p y}-\f{\p \hat{g}_*}{\p y}\|^2_{L^2(D)}
\lesssim \lambda.
$$
Similarly, we can show that
$$
\|\f{\p \hat{g}}{\p t}-\f{\p \hat{g}_*}{\p t}\|^2_{L^2(D)}
\lesssim \lambda.
$$
This completes the proof of (\ref{g-approximation1-1}) and hence the lemma.

\medskip

\textbf{Proof of Lemma \ref{lem-g-sum}}. Denote $D=(3\epsilon_1/4, t_1+ \epsilon_1/2)\times U(x_1)$ again.
We first show (\ref{gb-boundary1}). Since $x$ is restricted to $V(x_1)\subset \Gamma$, it suffices to show that
$$
\hat{g}^{-}(t, y, \lambda)+  \hat{g}(t, y, \lambda) =
O(\sqrt{\lambda}) \quad \mbox{in \,\,} H^1(D).
$$
But this is a direct consequence of Lemma \ref{lem-g-ghat1} and the fact that $\hat{g}_*^{-}=-\hat{g}_*$.

We now prove (\ref{gb-boundary2}). By direct calculate
\beas
\f{\p g}{\p \nu}(t, x)
&=&  \f{\p}{\p \nu}(a(t)e^{i\lambda\tau(t, x)})= i \lambda g\cdot \f{\p \tau}{\p \nu} \\
&=&   i \lambda g\cdot \langle \xi(t)+ M(t)(x-x(t), \nu(x) \rangle \\
&=&   i \lambda g \cdot \langle \xi(t_1), \nu(x_1) \rangle+  i \lambda g\cdot (\langle \xi(t), \nu(x) \rangle-  \langle \xi(t_1), \nu(x_1) \rangle)
      + i \lambda g \cdot \langle M(t)(x-x(t), \nu(x) \rangle.
\eeas
Note that in the coordinate $x= F(y)$,
\beas
|\langle M(t)(x-x(t)), \nu(x) \rangle| &=& O(|(t-t_1, y)|), \\
|\langle \xi(t), \nu(x) \rangle-  \langle \xi(t_1), \nu(x_1) \rangle|&=& O(|(t-t_1, y)|).
\eeas
It follows that
\beas
\|g\cdot (\langle \xi(t), \nu(x) \rangle-  \langle \xi(t_1), \nu(x_1) \rangle)\|_{L^2(D)}^2 & \lesssim &\f{1}{\lambda}, \\
\|g\cdot \langle M(t)(x-x(t), \nu(x) \rangle \|_{L^2(D)}^2 & \lesssim & \f{1}{\lambda}.
\eeas
Thus
$$
\f{\p g}{\p \nu}(t, x) = i \lambda g \cdot \langle \xi(t_1), \nu(x_1) \rangle + O(\sqrt{\lambda}).
$$
Similarly,
$$
\f{\p g^-}{\p \nu}(t, x) = i \lambda g^- \cdot \langle \xi^-(t_1), \nu(x_1) \rangle + O(\sqrt{\lambda}).
$$
Finally, using (\ref{gb-boundary1}) and the fact that $\langle \xi^-(t_1), \nu(x_1) \rangle= -\langle \xi(t_1), \nu(x_1) \rangle$, we conclude that (\ref{gb-boundary2}) holds. This completes the proof of the lemma.

\subsection{An estimate on Gaussian beam interactions} \label{appendix-beam-interaction}

\begin{lem} \label{inequality2c}
Assume that $M_1$ and $M_2$ are two symmetric positive definite matrices such that $ 0<c_0< M_1, M_2< c_1 $, and $N_1$ and
$N_2$ are two symmetric matrices such that $\|N_1\|, \|N_2\| \leq c_2$. Let
$\delta x, \delta \xi$  be two vectors in $\R^d$ and $\lambda \gg 1$. Then exists $c_3>0$ depending only on $c_0, c_1$ and $c_2$ such that
\begin{eqnarray*}
 |\int_{\R^d}  e^{i \lambda \cdot \langle \delta \xi, x \rangle -\lambda x^T(M_1+ i \cdot N_1)x- \lambda (x-\delta x)^T(M_2+ i \cdot N_2)(x-\delta x)}|  \lesssim \f{1}{\lambda^{\f{d}{2}}}e^{-c_3\lambda (|\delta x|^2+ |\delta \xi|^2)}.
\end{eqnarray*}
\end{lem}

{\bf{Proof.}} See Lemma 3.7 in \cite{BQ1}.


\begin{thebibliography} {99}
\bibitem{AS90}
G.\,Alessandrini and J.\,Sylvester.
Stability for a multidimensional inverse spectral theorem, Comm.PDE, 15(1990), 711-736.


\bibitem{AAB}
R.\,Alexandre, J.\,Akian and S.\,Bougacha.
Gaussian beams summation for the wave equation in a convex domain,
Commun. Math. Sci., 7(4)(2009), 973-1008.

\bibitem{A}
V.\,I.\,Arnold.
Singularity Theory, London Mathematical Society Lecture Note Series(53), Cambrige Univeristy Press, 1981.

\bibitem{AZV}
V.\,I.\,Arnold, S.\,M.\,Gusein-Zade and A.\,N.\,Varchenko.
Singularites of Differentiable Maps I, volume 82 of monographs in mathematics, Birkh$\ddot{a}$user, Boston, 1985.


\bibitem{BQ1}
G.\,Bao, J.\,Qian, L.\,Ying and H.\,Zhang.
A convergent multiscale Gaussian-beam parametrix for wave equations, Comm.PDE, 38(2013), 92-134.

\bibitem{BY09}
G.\,Bao and K.\,Yun.
On the stability of an inverse problem for the wave equation,
Inverse Problems, 25(4)(2009), R1-R7.



\bibitem{B87}
M.\,Belishev.
An approach to multidimensional inverse problems for the wave equation, Dokl.\,Akad.\,Nauk. SSR 297(1987), 524-527.

\bibitem{B97}
M.\,Belishev.
Boundary control in reconstruction of manifolds and metrics(the BC method), Inverse Problems, 13(1997), R1-R45.

\bibitem{B07}
M.\,Belishev.
Recent progress in the boundary control method,
Inverse Problems, 23(2007), R1-R67.

\bibitem{BK92}
M.\,Belishev and Y.\,V.\,Kurylev.
To the reconstruction of a Riemannian manifold via its spectral data(BC-method), Comm.P.D.E, 17(1992), 765-804.

\bibitem{BF11}
M.\,Bellassoued and D.\,D.\,S.\,Ferreira.
Stability estimates for the anisotropic wave equation from the Dirichlet-to-Neumann map,
Inverse Problems and Imaging, 5(4)(2011), 745-773.

\bibitem{GS-AJM-98}
A.\,Greenleaf and A.\,Seeger.
Fourier integral operators with cusp singularites,
American Journal of Mathematics, 120(5)(1998), 1077-1119.

\bibitem{H-III85}
L.\,H$\ddot{o}$mander.
The Analysis of Linear Partial Differential Operators, III and IV, volume 274 and 275 of Grundelehren der Mathematischen Wissenschaften, Spring-Verlag, Berlin, 1985.


\bibitem{Isakov06}
V.\,Isakov.
Inverse Problems for Partial Differential Equations, volume 127 of Applied Mathematical Sciences, second edition,
Springer Science+Business Media, Inc, New York, 2006.

\bibitem{KKL01}
A.\,Katchalov, Y.\,Kurylev and M.\,Lassas.
Inverse Boundary Spectral Problems, Chapman Hall/CRC, Boca Raton, 2001.

\bibitem{LLT86}
I.\,Lasiecka, J.\,L.\,Lions and R.\,Triggiani.
Non homogeneous boundary value problems for second order hyperbolic operators,
J. Math. Pures et Appl., 65(1986), 149-192.

\bibitem{Lee}
J.\,M.\,Lee.
Riemannian manifolds, an introduction to curvature, volume 176 of Graduate Texts in Mathematics, Spring-Verlag, New York, 1997.

\bibitem{Mi81}
R.\,Michel. Sur la rigidit\'e impos\'ee par la longueur des g\'eod\'esiques, Invent. Math. 65(1981), 71-83.

\bibitem{M-arxiv}
C.\,Montalto.
Stable determination of a simple metric, a covector field and a potential from the hyperbolic Dirichlet-to-Nuemann map,
accepted by Comm.PDE.


\bibitem{QY1}
J.\,Qian and L.\,Ying.
Fast multiscale Gaussian wavepacket transforms and multiscale Gaussian beams for the wave equation, SIAM. MMS, 8(2010).


\bibitem{PU06}
L.\,Pestov and G.\,Uhlmann.
The scattering relation and the Dirichlet-to-Neumann map,
Contemporary Math, 412(2006), 249-262.

\bibitem{R83}
J.\,Ralston.
Gassian beams and the propagation of singularities,
Studies in Partial Differential Equations, MAA Studies in Mathematics, Edited by W.\,Littman, 23(1983), 206-248.


\bibitem{R87}
V.\,G.\,Romanov.
Invere Problems of Mathematical Physics, VNU Science Press, Utreht, the Netherlands, 1987.

\bibitem{S99}
V.\,A.\,Sharafutdinov.
Ray Transform on Riemannian Manifolds, Lecture notes, University of Oulu, 1999.


\bibitem{SU98}
P.\,Stefanov and G.\,Uhlmann.
Stability estimates for the hyperbolic Dirichlet to Neumann map in anisotropic media,  J. Funct. Anal., 154(2)(1998), 330-358.



\bibitem{SU04Duke}
P.\,Stefanov and G.\,Uhlmann.
Stability estimates for the X-ray transform of tensor fields and boundary rigidity,
Duke Math. J., 123(2004), 445-467.


\bibitem{SU05JAMS}
P.\,Stefanov and G.\,Uhlmann.
Boundary rigidity and stability for generic simple metrics, Journal Amer. Math. Soc., 18(2005), 975-1003.

\bibitem{SU05IMRN}
P.\,Stefanov and G.\,Uhlmann.
Stable determination of generic simple metrics from the hyperbolic Dirichlet-to-Neumann map, IMRN, 17(2005), 1047-1061.

\bibitem{SU08JGA}
P.\,Stefanov and G.\,Uhlmann.
The X-Ray transform for a generic family of curves and weights,
J. Geom. Anal., 18(1)(2008), 81-97

\bibitem{SU08AJM}
P.\,Stefanov and G.\,Uhlmann.
Integral geometry of tensor fields on a class of non-simple Riemannian manifolds,
American J. of Math., 130 (2008), 239-268.



\bibitem{SU09JFA}
P.\,Stefanov and G.\,Uhlmann.
Linearizing non-linear inverse problems and its applications to inverse backscattering,
Journal Functional Analysis, 256 (2009), 2842-2866.

\bibitem{SU09JDG}
P.\,Stefanov and G.\,Uhlmann.
Local lens rigidity with incomplete data for a class of non-simple Riemannian manifolds,
Journal Diff. Geometry, 82 (2009), 383-409.

\bibitem{S90}
Z.\,Sun.
On continuous dependence for an inverse initial boundary value problem for the wave equation,
J. Math. Anal. Appl., 150(1990), 188-204.

\bibitem{U03}
G.\,Uhlmann.
The cauchy data and the scattering relation, IMA Publications, Springer-Verlag, 137(2003), 263-288.

\bibitem{SUAPDE}
P.\,Stefanov and G.\,Uhlmann.
The geodesic X-ray transform with fold caustics,
Anal. and PDE, 5(2012), 219-260.

\bibitem{thesis-zhang}
H.\,Zhang.
On the stability/sensitivity of recovering velocity fields from boundary measurements,
(PhD dissertation), Michigan State University, 2013.

\end{thebibliography}
\end{document}